\input amstex \documentstyle{amsppt} \nologo
\loadbold
\let\bk\boldkey

\hsize=5.75truein
\vsize=8.75truein 
\hcorrection{.25truein}
\loadeusm \let\scr\eusm
\loadeurm 
\font\Smc=cmcsc10 scaled \magstep1
\font\Rm=cmr17
\font\It=cmmi12
\define\cind{\text{\it c-\/\rm Ind}} 
\define\Ind{\text{\rm Ind}}
\define\Aut#1#2{\text{\rm Aut}_{#1}(#2)}
\define\End#1#2{\text{\rm End}_{#1}(#2)}
\define\Hom#1#2#3{\text{\rm Hom}_{#1}({#2},{#3})} 
\define\GL#1#2{\roman{GL}_{#1}(#2)}
\define\M#1#2{\roman M_{#1}(#2)}
\define\Gal#1#2{\text{\rm Gal\hskip.5pt}(#1/#2)}
\define\Ao#1#2{\scr A_{#1}(#2)} 
\define\Go#1#2{\scr G_{#1}(#2)} 
\define\upr#1#2{{}^{#1\!}{#2}} 
\define\N#1#2{\text{\rm N}_{#1/#2}} 
\define\sw{\text{\rm sw}}
\let\ge\geqslant
\let\le\leqslant
\let\ups\upsilon
\let\vD\varDelta 
\let\ve\varepsilon 
\let\eps\epsilon 
\let\vf\varphi
\let\vF\varPhi 
\let\vG\varGamma

\let\vp\varpi
\let\vP\varPi 
\let\vPs\varPsi 
 
\let\vs\varsigma 
\let\vt\vartheta
\let\vT\varTheta
 
\let\vX\varXi
\define\wP#1{\widehat{\scr P}_{#1}} 
\define\wW#1{\widehat{\scr W}_{#1}} 
\define\wG#1{\widehat{\text{\rm GL}}_{#1}} 
\topmatter 
\title\nofrills \Rm
Langlands parameters for epipelagic representations of GL$_{\text{\It n}}$ 
\endtitle 
\rightheadtext{Epipelagic Langlands parameters} 
\author 
Colin J. Bushnell and Guy Henniart 
\endauthor 
\leftheadtext{C.J. Bushnell and G. Henniart}
\affil 
King's College London and Universit\'e de Paris-Sud 
\endaffil 
\address 
King's College London, Department of Mathematics, Strand, London WC2R 2LS, UK. 
\endaddress
\email 
colin.bushnell\@kcl.ac.uk 
\endemail
\address 
Institut Universitaire de France et Universit\'e de Paris-Sud, Laboratoire de Math\'ematiques d'Orsay,
Orsay Cedex, F-91405; CNRS, Orsay cedex, F-91405. 
\endaddress 
\email 
Guy.Henniart\@math.u-psud.fr 
\endemail 
\date February 2013 \enddate 
\abstract 
Let $F$ be a non-Archimedean local field. An irreducible cuspidal representation of $\text{\rm GL}_n(F)$ is epipelagic if its Swan conductor equals $1$. We give a full and explicit description of the Langlands parameters of such representations. 
\endabstract 
\endtopmatter 
\document \nopagenumbers 
\baselineskip=14pt \parskip=4pt plus 1pt minus 1pt 
We examine an extreme class of irreducible cuspidal representations of the general linear group $G = \GL nF$ over a non-Archimedean local field $F$, those designated ``epipelagic'' by Reeder and Yu \cite{26}. A particular arithmetic interest of such representations was first identified by Gross and Reeder in \cite{18}, who dubbed them ``simple cuspidals''. The authors of \cite{18} asked us whether anything could be said about the Langlands parameters. This paper is the response. 
\par 
The context of \cite{26} is quite general but, for the group $\GL nF$, an exercise in techniques going back to \cite{1} shows that a cuspidal representation $\pi$ is epipelagic if and only if the exponential Swan conductor $\roman{Sw}(\pi)$ is equal to $1$. We start from that point. 
\subhead 
1 
\endsubhead 
If $n\ge1$ is an integer, let $\Ao nF$ be the set of equivalence classes of irreducible, cuspidal (complex) representations of $G = \GL nF$. Let $\Go nF$ be the set of equivalence classes of irreducible, smooth, $n$-dimensional representations of the Weil group $\scr W_F$ of some separable algebraic closure $\bar F/F$ of $F$. We denote the Langlands correspondence $\Ao nF \to \Go nF$ by $\pi\mapsto \upr L\pi$. Let $\upr 1{\Ao nF}$ be the set of $\pi\in \Ao nF$ with $\roman{Sw}(\pi) = 1$ and define $\upr 1{\Go nF}$ in the same way. The Langlands correspondence preserves Swan conductors and so maps $\upr 1{\Ao nF}$ bijectively to $\upr 1{\Go nF}$. 
\par
It is easy to describe the elements of $\upr 1{\Ao nF}$ in terms of the standard model of \cite{14}. A representation $\pi \in \Ao nF$, not of level zero, contains a simple character $\theta$. If we choose a character $\psi_F$ of $F$ satisfying $c(\psi_F) = -1$ (see \S1 for this notation), then $\theta$ is attached to a simple stratum $[\frak a,l,0,\alpha]$ in the matrix algebra $A = \M nF$. The condition $\roman{Sw}(\pi) = 1$ translates into $l$ being $1$, the hereditary order $\frak a$ being minimal and the field extension $F[\alpha]/F$ being totally ramified of degree $n$. Fixing such a stratum $[\frak a,1,0,\alpha]$, one may list the elements $\pi$ of $\upr 1{\Ao nF}$ to which it is attached. This yields a classification of the elements of $\upr1{\Ao nF}$ in terms of properties of local constants. Passing across the Langlands correspondence, one gets a parallel classification of the elements of $\upr 1{\Go nF}$. 
\subhead 
2 
\endsubhead 
On the Galois side, that classification exhibits a lamentable opacity. We set ourselves the task of describing explicitly the representation $\sigma = \upr L\pi$, $\pi \in \upr 1{\Ao nF}$, in terms of the structure of $\pi$. 
\par
There are not many cases where such descriptions have been seriously attempted. The essentially tame case is exhaustively worked out in \cite{7}, \cite{8}, \cite{10}. The wildly ramified case has only been considered in prime dimension \cite{9}, \cite{20}, Kutzko \cite{24}, M\oe glin \cite{25}, with the aim of verifying the Langlands conjecture, as it was at the time. The dissection of the representations was carried no further than necessary for that purpose. A thorough analysis, even of a special case, seems overdue. The epipelagic representations exhibit sufficient complexity to make the exercise worthwhile and can be used to generate further families of examples. 
\subhead 
3 
\endsubhead 
Let $\sigma\in \upr 1{\Go nF}$ and $\pi\in \upr1{\Ao nF}$ satisfy $\upr L\pi = \sigma$. Let $p$ be the residual characteristic of $F$, and write $n = ep^r$, for integers $e$, $r$ with $e$ not divisible by $p$. There is then a totally ramified extension $K/F$ of degree $e$, and a representation $\tau\in \upr1{\Go{p^r}F}$, such that $\tau$ induces $\sigma$. The pair $(K/F,\tau)$ is uniquely determined up to conjugation in $\scr W_F$. The theory of tame lifting enables us to specify $K$ in terms of $\pi$, and also the representation $\rho\in \upr 1{\Ao{p^r}K}$ such that $\tau  = \upr L\rho$. 
\subhead 
4 
\endsubhead 
We therefore concentrate on the case $\upr L\pi = \sigma \in \upr1{\Go{p^r}F}$. Here, a conductor estimate shows that $\sigma$ is {\it primitive,} and so covered by Koch's seminal work \cite{23}. 
\par
Let $\bar\sigma$ be the {\it projective\/} representation of $\scr W_F$ determined by $\sigma$, and define a finite Galois extension $E/F$ by $\scr W_E = \roman{Ker}\,\bar\sigma$. The restriction of $\sigma$ to $\scr W_E$ is therefore a multiple of a character $\xi_\sigma$ of $\scr W_E$. Let $T/F$ be the maximal tamely ramified sub-extension of $E/F$. Following \cite{23}, the restriction $\sigma_T$ of $\sigma$ to $\scr W_T$ is irreducible. The Galois group of $E/T$ is elementary abelian of order $p^{2r}$, and $\sigma_T$ is the unique irreducible representation of $\scr W_T$ containing $\xi_\sigma$. The classification theory of \cite{23} then rests on the fact that $\Gal ET$ provides a symplectic representation of $\Gal TF$ over the field $\Bbb F_p$ of $p$ elements. 
\par 
There is a helpful element of structure here. If $K/F$ is a finite tame extension, let $\sigma_K$ denote the (irreducible) restriction of $\sigma$ to $\scr W_K$. Let $D(\sigma_K)$ be the group of characters $\chi$ of $\scr W_K$ for which $\chi\otimes \sigma_K \cong \sigma_K$. One then has $|D(\sigma_K)| \le p^{2r}$, with equality if and only if $K$ contains $T$. Moreover, $D(\sigma_T)$ is the dual of the Galois group $\Gal ET = \scr W_T/\scr W_E$. 
\par 
We can make no further progress working only with Galois representations and have to switch to the other side. If $K/F$ is a finite tame extension, we may define $\pi_K\in \Ao{p^r}K$ by the relation $\upr L\pi_K = \sigma_K$. Using results from \cite{4}, \cite{5} on tame lifting, we can describe $\pi_K$ explicitly, in terms of $\pi$. We have to find the least tame extension $T/F$ for which there are $p^{2r}$ characters $\phi$ of $T^\times$ satisfying $\phi\pi_T \cong \pi_T$. An elaborate calculation is required. The outcome is an explicit polynomial, with splitting field $T$. The roots of the polynomial define the non-trivial elements of $D(\sigma_T)$, as characters of $T^\times$. Via class field theory, they determine the abelian extension $E/T$. 
\par 
The structure of the argument is worthy of comment. The initial steps, from the Galois theory of local fields, are fairly straightforward. We only use the early parts of \cite{23}, not going much beyond the group-theoretic arguments of Rigby \cite{27}. The main calculations, leading to the polynomial determining the field $T$, are quite involved and use a lot of machinery from \cite{14}, but they are essentially elementary in nature. The depth lies in the facility with which we can move from side to side, via the Langlands correspondence, without loss of explicit detail. This is absolutely reliant on a substantial portion of the theory of tame lifting, as developed in \cite{4}, \cite{6} and \cite{12}.  
\subhead 
5 
\endsubhead  
Having identified the fields $T/F$ and $E/T$, it remains to understand the ``p-central character'' $\xi_\sigma$. Starting from the description of $T$ and $D(\sigma_T)$, classical methods of local field theory and a local constant calculation yield an expression for $\xi_\sigma$ on the unit group $U^1_E$. We shall see that the datum $(\xi_\sigma|_{U^1_E}, \det\sigma)$ determines $\sigma$ up to tensoring with an unramified character of order dividing $p^r$ (6.2 below). Since  $\roman{Sw}(\sigma)$ is not divisible by $p$, this uncertainty is removed by a simple local constant relation. 
\par 
It is possible to specify $\xi_\sigma$ completely by calculating one further special value, but the details are voluminous. Nothing much is gained, so we have omitted them. 
\head \Smc 
1. Preliminaries and notation
\endhead 
We set down our standard notation and recall some basic facts. 
\subhead 
1.1 
\endsubhead 
Throughout, $F$ is a non-Archimedean local field with discrete valuation ring $\frak o_F$. The maximal ideal of $\frak o_F$ is $\frak p_F$, the residue field is $\Bbbk_F = \frak o_F/\frak p_F$, $q = |\Bbbk_F|$, and $\ups_F:F^\times \twoheadrightarrow \Bbb Z$ is the normalized additive valuation. The characteristic of $\Bbbk_F$ is $p$. The unit group is $U_F = \frak o_F^\times$ and $U^k_F = 1{+}\frak p_F^k$, $k\ge1$. Also, $\mu_F$ is the group of roots of unity in $F$, of order relatively prime to $p$. 
\par 
If $\frak a$ is a hereditary $\frak o_F$-order in some matrix algebra $\M nF$ and $\frak p_\frak a$ is the Jacobson radical of $\frak a$, then  $U_\frak a = \frak a^\times$ and $U^k_\frak a = 1{+}\frak p_\frak a^k$, $k\ge1$. 
\par 
Let $\bar F/F$ be a separable algebraic closure of $F$. Let $\scr W_F$ be the Weil group of $\bar F/F$ and $\scr P_F$ the wild inertia subgroup of $\scr W_F$. We identify the group of smooth characters of $\scr W_F$ with that of $F^\times$ via class field theory, switching between the two viewpoints as convenient.
\par 
Let $\wW F$ be the set of equivalence classes of irreducible smooth representations of $\scr W_F$. For an integer $n\ge1$, let $\Go nF$ be the set of $\sigma\in \wW F$ of dimension $n$. Let $\Ao nF$ be the set of equivalence classes of irreducible {\it cuspidal\/} representations of $\GL nF$. The Langlands correspondence gives a bijection $\Ao nF \to \Go nF$ which we denote by $\pi\mapsto \upr L\pi$. 
\subhead 
1.2 
\endsubhead 
Let $\psi$ be a smooth character of $F$, $\psi\neq1$. Define $c(\psi)$ to be the greatest integer $k$ such that $\frak p_F^{-k} \subset \roman{Ker}\,\psi$. We shall mostly work with the case $c(\psi) = -1$, where $\psi$ is trivial on $\frak p_F$ but not on $\frak o_F$. 
\par 
Let $n\ge1$ be an integer, let $\pi\in \Ao nF$. The {\it Godement-Jacquet local constant\/} $\ve(\pi,s,\psi)$ of $\pi$, relative to the character $\psi$ and a complex variable $s$, takes the form (Godement-Jacquet \cite{17},  Jacquet \cite{21}) 
$$
\ve(\pi,s,\psi) = q^{-s\roman f(\pi,\psi)}\ve(\pi,0,\psi), 
$$ 
for an integer $\roman f(\pi,\psi)$ such that 
$$ 
\roman f(\pi,\psi) = \roman a(\pi)+nc(\psi), 
$$ 
where the {\it Artin conductor\/} $\roman a(\pi)$ is independent of $\psi$. Set 
$$
\sw(\pi) = \roman a(\pi)-n. 
$$
Thus $\sw(\pi) = \roman f(\pi,\psi)$ when $c(\psi) = -1$. The integer $\sw(\pi)$ is the same as the Swan conductor $\roman{Sw}(\pi)$ {\it except\/} when $n=1$ and $\pi$ is an unramified character of $F^\times$. In that case, $\sw(\pi) = -1$ and $\roman{Sw}(\pi) = 0$. 
\par 
We denote by $\upr1{\Ao nF}$ the set of $\pi\in \Ao nF$ for which $\sw(\pi) = 1$.
\subhead 
1.3 
\endsubhead 
Let $\sigma\in \Go nF$. The {\it Langlands-Deligne local constant\/} $\ve(\sigma,s,\psi)$ takes the form \cite{9}, Tate \cite{30}, 
$$
\ve(\sigma,s,\psi) = q^{-s\roman f(\sigma,\psi)}\ve(\sigma,0,\psi), 
$$ 
with $\roman f(\sigma,\psi) = \roman a(\sigma)+nc(\psi)$, the Artin conductor $\roman a(\sigma)$ being independent of $\psi$. Likewise set 
$$
\sw(\sigma) = \roman a(\sigma)-n, \quad n = \dim\sigma.  
$$ 
This relates to the Swan conductor $\roman{Sw}(\sigma)$ as before. Let $\upr1{\Go nF}$ be the set of $\sigma\in \Go nF$ with $\sw(\sigma) = 1$. The Langlands correspondence satisfies $\sw(\upr L\pi) = \sw(\pi)$, $\pi \in \Ao nF$, so it maps $\upr 1{\Ao nF}$ bijectively to $\upr 1{\Go nF}$. 
\par 
We make frequent use of the following fact. 
\proclaim{Induction Formula} 
Let $E/F$ be a finite separable extension, set $e=e(E|F)$, $f=f(E|F)$, and let $\frak D_{E/F} = \frak p_E^d$ be the relative different of $E/F$. If $\tau\in \wW E$ has dimension $m$ and $\sigma = \Ind_{E/F}\,\tau$, then 
$$
\sw(\sigma) = \big(\sw(\tau)+m(1{-}e{+}d)\big)f. 
\tag 1.3.1
$$
\endproclaim 
\demo{Proof} 
If $\psi$ is a smooth character of $F$, $\psi\neq1$, and $\psi_E = \psi\circ \roman{Tr}_{E/F}$, then 
$$ 
c(\psi_E) = d +  ec(\psi). 
\tag 1.3.2 
$$ 
According to \cite{9} 30.4 Corollary, the quotient $\ve(\sigma,s,\psi)/\ve(\tau,s,\psi_E)$ is independent of $s$. The result follows from a simple computation. \qed 
\enddemo  
\head\Smc 
2. Classification
\endhead 
Let $n\ge1$ be an integer. In this section, we classify the elements of $\upr 1{\Ao nF}$ using the method of \cite{14}. We translate this into a classification in terms of properties of local constants. Via the Langlands correspondence, this gives a parallel classification of the elements of $\upr 1{\Go nF}$. 
\par 
Since we work in the context of \cite{14}, we need to choose a smooth character $\psi_F$ of $F$ with $c(\psi_F) = -1$. Thus, in the language of \cite{14}, the character $\psi_F$ ``has level one''. Throughout, we write $G = \GL nF$ and $A = \M nF$. 
\subhead 
2.1 
\endsubhead 
Let $\pi\in \Ao nF$. Thus $\pi$ contains a {\it simple character\/} $\theta$ in $G$, in the sense of \cite{14}. This simple character is uniquely determined up to $G$-conjugation. It is trivial if and only if $\sw(\pi) = 0$. Assuming $\theta$ to be non-trivial, it is attached via $\psi_F$ to a simple stratum $[\frak a,l,0,\beta]$ in $A$: in the notation of \cite{14}, $\theta\in \scr C(\frak a,\beta,\psi_F)$. The stratum $[\frak a,l,0,\beta]$ is {\it m-simple,} in that $\frak a$ is maximal among hereditary $\frak o_F$-orders in $A$, stable under conjugation by $F[\beta]^\times$. Correspondingly, one says that $\theta$ is m-simple. 
\par 
If $e_\frak a$ denotes the $F$-period of the hereditary $\frak o_F$-order $\frak a$, then 
$$
\sw(\pi) = ln/e_\frak a \ge 1, 
$$
\cite{5} 6.1 Lemma 2. Since the stratum $[\frak a,l,0,\beta]$ is m-simple, the period $e_\frak a$ equals the ramification index $e(F[\beta]|F)$ of the field extension $F[\beta]/F$. We encapsulate these remarks:  
\proclaim{Lemma 1} 
Let $\pi \in \Ao nF$. The following conditions are equivalent: 
\roster 
\item $\sw(\pi) = 1$; 
\item $\pi$ contains an m-simple character $\theta \in \scr C(\frak a,\alpha,\psi_F)$, where $[\frak a,1,0,\alpha]$ is an m-simple stratum in $A$ such that $e_\frak a = n$. 
\endroster 
When these conditions are satisfied, the field extension $F[\alpha]/F$ is totally ramified of degree $n$. 
\endproclaim 
We write down the various groups and the simple characters attached to the simple stratum $[\frak a,1,0,\alpha]$ in the machinery of \cite{14} Chapter 3. Let $\roman{tr}_A:A\to F$ be the matrix trace and write $\psi_A = \psi_F\circ \roman{tr}_A$. 
\proclaim{Lemma 2} 
Let $[\frak a,1,0,\alpha]$ be an m-simple stratum in $A$ with $e_\frak a = n$. 
\roster 
\item 
The associated groups are 
$$
H^1(\alpha,\frak a) = J^1(\alpha,\frak a) = U^1_\frak a. 
$$
The set $\scr C(\frak a,\alpha,\psi_F)$ of simple characters has only the one element 
$$
\theta_\alpha:1{+}x \longmapsto \psi_A(\alpha x), \quad 1{+}x\in U^1_\frak a. 
$$
The $G$-normalizer of the character $\theta_\alpha$ is $\bk J_\alpha = F[\alpha]^\times U^1_\frak a$. 
\item 
Let $\vPs$ be a character of $\bk J_\alpha$ such that $\vPs|_{U^1_\frak a} = \theta_\alpha$. The $G$-representation 
$$ 
\pi_\vPs = \cind_{\bk J_\alpha}^G\vPs 
$$
is irreducible and cuspidal. The map $\vPs \mapsto \pi_\vPs$ is a bijection between the set of characters of $\bk J_\alpha$ extending $\theta_\alpha$ and the set of $\pi\in \Ao nF$ containing $\theta_\alpha$. 
\endroster 
\endproclaim 
\demo{Proof} 
Part (1) follows on working through the definitions in Chapter 3 of \cite{14}. Part (2) is an instance of \cite{14} (8.4.1) or \cite{12} 4.1. \qed 
\enddemo 
In part (1) of the lemma, the character $\theta_\alpha$ determines only the coset $\alpha{+}\frak a = \alpha U^1_\frak a$. 
In the situation of part (2), it will be easier to say that $\pi$ {\it contains the m-simple stratum} $[\frak a,1,0,\alpha]$ (relative to the character $\psi_F$), rather than that it contains $\theta_\alpha$. 
\proclaim{Proposition} 
\roster 
\item If $\pi \in \upr 1{\Ao nF}$, $n>1$, the central character $\omega_\pi$ of $\pi$ is tamely ramified. 
\item Let $\omega$ be a tamely ramified character of $F^\times$ and $[\frak a,1,0,\alpha]$ an m-simple stratum in $A$, with $e_\frak a = n$. There exist exactly $n$ representations $\pi \in \upr 1{\Ao nF}$ containing $[\frak a,1,0,\alpha]$ and having central character $\omega_\pi = \omega$. 
\endroster 
\endproclaim 
\demo{Proof} 
Let $\pi \in \upr 1{\Ao nF}$ contain the m-simple stratum $[\frak a,1,0,\alpha]$. We have $U^1_F = F^\times \cap U^n_\frak a \subset \roman{Ker}\,\theta_\alpha$. Therefore $\omega_\pi$ is trivial on $U^1_F$, as required for (1). The group $\bk J_\alpha$ contains $F^\times U^1_\alpha$ with index $n$, so (2) follows from Lemma 2. \qed 
\enddemo 
\subhead 
2.2 
\endsubhead  
If $\pi \in \Ao nF$ and if $\chi$ is a character of $F^\times$, then $\chi\pi$ will denote the representation $g\mapsto \chi(\det g)\pi(g)$. In particular, $\chi\pi \in \Ao nF$. 
\par 
We compute Godement-Jacquet local constants. The method originates in \cite{3}, \cite{1}, but the summary in \cite{5} 6.1 is closer to our current notation. 
\proclaim{Lemma} 
If $\pi = \cind_{\bk J_\alpha}^G\vPs_\pi$, as in \rom{2.1 Lemma 2,} then 
\roster 
\item $\ve(\pi,\tfrac12,\psi_F) = \vPs_\pi(\alpha)^{-1}\,\psi_A(\alpha)$ and,  
\item if $\chi$ is a tamely ramified character of $F^\times$, then 
$$ 
\align 
\sw(\chi\pi) &= \sw(\pi) = 1, \\ 
\ve(\chi\pi,s,\psi_F) &= \chi(\det\,\alpha)^{-1}\,\ve(\pi,s,\psi_F). 
\endalign 
$$ 
\endroster 
\endproclaim 
We may now strengthen the description of 2.1. 
\proclaim{Proposition} 
For $i=1,2$, let $\pi_i\in \upr 1{\Ao nF}$, and let $\omega_i$ be the central character of $\pi_i$. Let $\pi_i$ contain the m-simple stratum $[\frak a,1,0,\alpha_i]$. The representations $\pi_1$, $\pi_2$ are equivalent if and only if the following conditions are satisfied: 
\roster 
\item"\rm (a)" $\det\,\alpha_1 \equiv \det\alpha_2 \pmod{U^1_F}$, 
\item"\rm (b)" $\omega_1 = \omega_2$, and 
\item"\rm (c)" $\ve(\pi_1,\frac12,\psi_F) = \ve(\pi_2,\frac12,\psi_F)$. 
\endroster 
\endproclaim 
\demo{Proof} 
Let $\frak p_\frak a$ be the Jacobson radical of $\frak a$. If the $\pi_i$ are equivalent, the m-simple characters $\theta_{\alpha_i}$ are conjugate in $G$ or, equivalently, the cosets $\alpha_i U^1_\frak a$ are $U_\frak a$-conjugate. To elucidate this condition, we use a matrix normal form calculation. We may assume that $\frak a$ is the ring of matrices, with entries in $\frak o_F$, which become upper triangular when reduced modulo $\frak p_F$. The ideal $\frak p_\frak a$ then consists of the matrices which are strictly upper triangular modulo $\frak p_F$. We have $\alpha_i^{-1}\frak a = \frak p_\frak a$ for both values of $i$, this being part of the definition of a simple stratum. A matrix manipulation shows that $\alpha_i^{-1}$ is $U_\frak a$-conjugate to a matrix $a_iu_i$, where 
\roster 
\item $u_i\in U^1_\frak a$; 
\item $\big(a_i\big)_{j,j+1} = 1$, for $1\le j<n$; 
\item $\big(a_i\big)_{n,1}$ is a prime element $\vp_i$ of $F$; 
\item all other entries of $a_i$ are zero. 
\endroster 
The cosets $a_iU^1_\frak a$ are then $U_\frak a$-conjugate if and only if $\vp_1\equiv \vp_2 \pmod{U^1_F}$, and this is equivalent to condition (a) of the proposition. 
\par
We are thus reduced to the case $\alpha_1 = \alpha_2 = \alpha$, say. Set $E = F[\alpha]$. There are characters $\vPs_i$ of $\bk J_\alpha = E^\times U^1_\frak a$ such that $\vPs_i|_{U^1_\frak a} =\theta_\alpha$ and $\pi_i \cong \cind_{\bk J_\alpha}^G \vPs_i$, $i=1,2$. Moreover, $\pi_1 \cong \pi_2$ if and only if $\vPs_1 = \vPs_2$. Condition (b) of the proposition is equivalent to $\vPs_1|_{F^\times} = \vPs_2|_{F^\times}$ (2.1 Lemma 2). Using the last lemma, condition (c) of the proposition is now equivalent to $\vPs_1(\alpha) = \vPs_2(\alpha)$. Since $\alpha$ generates the finite cyclic group $\bk J_\alpha/F^\times U^1_\frak a$, the result follows. \qed 
\enddemo 
\remark{Remark} 
Let $\omega$ be a tamely ramified character of $F^\times$. Taking account of 2.1 Proposition, we see that the set $\upr 1{\Ao nF}$ has exactly $n(q{-}1)$ elements $\pi$ such that $\omega_\pi = \omega$. 
\endremark 
\subhead 
2.3 
\endsubhead 
We use the Langlands correspondence to translate the classification from 2.2 Proposition in terms of representations of $\scr W_F$. 
\proclaim{Lemma} 
Let $\sigma\in \wW F$ and suppose that $\sw(\sigma) \ge1$. There exists $\gamma_\sigma\in F^\times$ such that 
$$
\ve(\chi\otimes\sigma,s,\psi_F) = \chi(\gamma_\sigma)^{-1}\,\ve(\sigma,s,\psi_F), 
$$
for any tamely ramified character $\chi$ of $\scr W_F$. This property determines the coset $\gamma_\sigma U^1_F$ uniquely. 
\endproclaim 
\demo{Proof} 
The lemma is a case of the main result of \cite{16}. \qed 
\enddemo 
\proclaim{Proposition} 
Let $\sigma\in \upr 1{\Go nF}$. Define $\pi\in \upr 1{\Ao nF}$ by $\upr L\pi = \sigma$. If $[\frak a,1,0,\alpha]$ is an m-simple stratum contained in $\pi$, then 
$$ 
\align  
\gamma_\sigma &\equiv \det\alpha \pmod{U^1_F}, \\ \det\sigma &= \omega_\pi, \\ 
\ve(\sigma,s,\psi_F) &= \ve(\pi,s,\psi_F). 
\endalign 
$$ 
Moreover, $\pi$ is the only element of $\Ao nF$ with these properties. 
\endproclaim 
\demo{Proof} 
The Langlands correspondence $\Ao nF \to \Go nF$ takes the central character to the determinant, it preserves the local constant and respects twisting with characters. The result therefore follows from the lemma and 2.2 Lemma. \qed 
\enddemo 
Combining the proposition with 2.2 Proposition,we get: 
\proclaim{Corollary} 
Let $\sigma_1,\sigma_2 \in \upr1{\Go nF}$. The representations $\sigma_1$, $\sigma_2$ are equivalent if and only if the following conditions hold: 
\roster 
\item $\gamma_{\sigma_1} \equiv \gamma_{\sigma_2} \pmod{U^1_F}$,  
\item $\det\sigma_1 = \det\sigma_2$, and 
\item $\ve(\sigma_1,s,\psi_F) = \ve(\sigma_2,s,\psi_F)$. 
\endroster 
\endproclaim 
\remark{Reflection} 
It is not difficult to reconstruct $\pi \in \upr 1{\Ao nF}$ from properties of local constants. In particular, the function $\chi\mapsto \ve(\chi\pi,\frac12,\psi)$, with $\chi$ ranging over the tame characters of $F^\times$, reveals completely the simple character contained in $\pi$. For $\sigma\in \upr1{\Go nF}$, the corresponding property gives the coset $\gamma_\sigma U^1_F$ and that determines the simple character in $\pi$, $\upr L\pi = \sigma$. According to the Ramification Theorem \cite{6} 8.2 or \cite{12} 6.1, it must also determine the restriction of $\sigma$ to the wild inertia subgroup of $\scr W_F$. The process by which it does this is the central concern of the paper. 
\endremark 
\head\Smc 
3. Primitivity 
\endhead 
While the results of \S2 classify the elements of $\upr1{\Go nF}$, they do not describe them at all effectively. We investigate further. 
\subhead 
3.1 
\endsubhead 
Any $\sigma\in \upr 1{\Go nF}$ is {\it totally ramified,} in the following sense. 
\proclaim{Lemma} 
Let $\sigma\in \upr 1{\Go nF}$ and let $\chi$ be an unramified character of $\scr W_F$. If the representations $\chi\otimes\sigma$, $\sigma$ are equivalent then $\chi = 1$. 
\endproclaim 
\demo{Proof} 
By 2.3 Proposition, $\ups_F(\gamma_\sigma) = -1$ in this case. If the unramified character $\chi$ is non-trivial, then  $\ve(\chi\otimes\sigma,s,\psi_F) = \chi(\gamma_\sigma)^{-1}\ve(\sigma,s,\psi_F) \neq \ve(\sigma,s,\psi_F)$, whence $\chi\otimes\sigma \not\cong \sigma$. \qed 
\enddemo 
Equivalently, a representation $\sigma\in \upr 1{\Go nF}$ restricts irreducibly to the inertia subgroup of $\scr W_F$. 
\proclaim{Proposition 1} 
Let $n = ep^r$, for integers $e$, $r$ such that $p$ does not divide $e$. Let $\sigma\in \upr 1{\Go nF}$. 
\roster 
\item 
There exists a totally ramified extension $K/F$, of degree $e$, and a representation $\tau\in \wW K$ such that $\sigma \cong 
\Ind_{K/F}\,\tau$. This relation determines the pair $(K/F,\tau)$ uniquely up to $\scr W_F$-conjugation. 
\item 
The representation $\tau$ satisfies $\sw(\tau) = 1$. 
\endroster 
\endproclaim 
\demo{Proof} 
Part (1) is an instance of \cite{6} 8.6 Proposition. Part (2) follows directly from the Induction Formula (1.3.1). \qed 
\enddemo 
We describe the representation $\tau$ in the manner of \S2. Let $\sigma = \upr L\pi$, $\pi\in \upr1{\Ao 
nF}$. As in 2.1, $\pi$ contains the simple character $\theta_\alpha\in \scr C(\frak a,\alpha,\psi_F)$, for a 
simple stratum $[\frak a,1,0,\alpha]$ in $A = \M nF$. 
\par
Let $T/F$ be the maximal tamely ramified sub-extension of $F[\alpha]/F$. Let $B \cong \M{p^r}T$ be the 
$A$-centralizer of $T$ and set $\frak b = \frak a\cap B$. Write $\psi_T = \psi_F\circ \roman{Tr}_{T/F}$. 
The quadruple $[\frak b,1,0,\alpha]$ is then an m-simple stratum in $B$ and the set $\scr C(\frak 
b,\alpha,\psi_T)$ has only one element $\theta_\alpha^T$, as in 2.1 Lemma 2. 
\par
We take $(K/F,\tau)$ as in Proposition 1. 
\proclaim{Proposition 2} 
\roster 
\item 
The fields $T$, $K$ are $F$-isomorphic. 
\item 
Let $\rho \in \upr 1{\Ao{p^r}K}$ satisfy $\upr L\rho = \tau$. There is an $F$-isomorphism $f:K\to T$ such 
that the representation $f_*\rho\in \upr 1{\Ao{p^r}T}$ contains the simple character $\theta_\alpha^T$. 
\endroster 
\endproclaim 
\demo{Proof} 
Part (1) follows from the Tame Parameter Theorem of \cite{12} 6.3 and the second from \cite{12} 6.2. 
\qed 
\enddemo 
To save notation, we identify $T$ with $K$ via the isomorphism $f$, and let $\det_B:B^\times \to K^\times$ be the determinant map. 
\proclaim{Corollary} 
Let $1_K$ be the trivial character of\/ $\scr W_K$. Set $R_{K/F} = \Ind_{K/F}\,1_K$ and $\delta_{K/F} = \det R_{K/F}$. The 
representation $\tau$ satisfies: 
\roster 
\item $\gamma_\tau \equiv \det_B\alpha \pmod{U^1_K}$; 
\item $\det\tau|_{F^\times} = \delta_{K/F}^{-p^r}\,\det_\sigma$; 
\item $\ve(\sigma,s,\psi_F)/\ve(\tau,s,\psi_K) = \big(\ve(R_{K/F},s,\psi_F)/\ve(1_K,s,\psi_K)\big)^{p^r}$.  
\endroster 
These relations determine $\tau$ uniquely. 
\endproclaim 
\demo{Proof}  
Part (1) is implied by Proposition 2 and 2.3 Proposition. For parts (2) and (3), see for instance \cite{9} 29.2, (29.4.1). Uniqueness follows from 2.3 Corollary. \qed 
\enddemo 
The values of the Langlands constant $\ve(R_{K/F},s,\psi_F)/\ve(1_K,s,\psi_K)$ are tabulated in, for example, \cite{2} 10.1. The corollary allows one to specify the representation $\rho$, such that $\upr L\rho = \tau$, in the manner of 2.2. 
\subhead 
3.2 
\endsubhead 
The results of 3.1 reduce us to the case where $\dim\sigma$ is a power of $p$. In particular, $\sigma$ is {\it totally wildly ramified,} in that it remains irreducible on restriction to the wild inertia subgroup $\scr P_F$ of $\scr W_F$. 
\par 
Recall that a representation $\tau\in \wW F$ is {\it imprimitive\/} if there exists a finite separable extension $K/F$, $K\neq F$, and a representation $\rho\in \wW K$ such that $\tau \cong \Ind_{K/F}\,\rho$. Thus $\tau$ is called {\it primitive\/} if it is not imprimitive. 
\proclaim{Proposition} 
If $\sigma\in \upr 1{\Go{p^r}F}$, for an integer $r\ge1$, then $\sigma$ is primitive. 
\endproclaim 
\demo{Proof} 
Suppose, for a contradiction, that $\sigma \cong \Ind_{K/F}\,\tau$, for a finite extension $K/F$, $K\neq F$, and $\tau\in \wW K$. Thus $[K{:}F] = p^a$ and $\dim\tau = p^b$, for integers $a,b$ such that $a{+}b = r$. By (1.3.1),  $\sw(\sigma) = 1$ is divisible by $f(K|F)$, so $K/F$ must be totally ramified. Since $[K{:}F]$ is a power of $p$, the extension $K/F$ is totally wildly ramified. If $\frak D_{K/F} = \frak p_K^{d(K|F)}$ is the different of $K/F$, this means $d(K|F) \ge [K{:}F] = p^a$. Therefore  
$$
1 = \sw(\sigma) \ge \sw(\tau) + p^b(1{-}p^a{+}p^a), 
$$
or 
$$
\sw(\tau) \le 1{-}p^b \le 0. 
$$ 
The only possibility here is $p^b = \dim\tau = 1$ and $\tau$ tamely ramified. The character $\tau$ of $K^\times$ is then of the form $\xi\circ\N KF$, for a tamely ramified character $\xi$ of $F^\times$. The character $\xi$ must occur as a component of $\Ind_{K/F}\,\tau$,  which is therefore reducible. This contradicts our hypothesis. \qed 
\enddemo 
\head\Smc 
4. Structure of primitive representations 
\endhead 
Following 3.2 Proposition, a representation $\sigma\in \upr 1{\Go{p^r}F}$, $r\ge1$, is necessarily primitive. Primitive representations are described in \cite{23}. We now see how $\sigma$ fits into the scheme of \cite{23}, as a basis for the detailed analysis of the next section. 
\subhead 
4.1 
\endsubhead 
We make some remarks related to Clifford theory. For $\tau\in \wW F$, let $D(\tau)$ denote the group of characters $\chi$ of $\scr W_F$ such that $\chi\otimes \tau\cong \tau$. 
\proclaim{Lemma} 
Let $\tau\in \wW F$ have dimension $n$. The abelian group $D(\tau)$ has exponent dividing $n$ and order at most $n^2$. 
\endproclaim 
\demo{Proof} 
Let $\chi\in D(\tau)$. The relation $\chi\otimes\sigma \cong \sigma$ implies $\det\sigma = \det(\chi\otimes\sigma) = \chi^n\det\sigma$, whence $\chi^n = 1$. Let $\check\tau$ denote the contragredient of $\tau$, and let $\chi$ be a character of $\scr W_F$. Frobenius Reciprocity yields 
$$ 
\Hom{\scr W_F}{\chi\otimes\tau}\tau \cong \Hom{\scr W_F}\chi{\check\tau\otimes\tau}. 
$$ 
We deduce that $\chi\in D(\tau)$ if and only if $\chi$ is a component of $\check\tau\otimes\tau$. This proves the second assertion. \qed 
\enddemo 
\proclaim{Proposition} 
Let $\tau\in \wW F$ be totally wildly ramified. 
\roster 
\item 
If $\chi \in D(\tau)$ is tamely ramified, then $\chi = 1$. 
\item 
If $K/F$ is a finite tame extension, the canonical map 
$$ 
\align 
D(\tau) &\longrightarrow D(\tau_K), \\ \chi&\longmapsto \chi_K = \chi\circ\N KF, 
\endalign 
$$ 
is injective. If $K/F$ is Galois, the image of $D(\tau)$ is the set of $\Gal KF$-fixed points in $D(\tau_K)$.
\endroster 
\endproclaim 
\demo{Proof} 
Let $\chi\in D(\tau)$ be tamely ramified, $\chi\neq 1$. The kernel of $\chi$ is then $\scr W_K$, for a cyclic tame extension $K/F$ and, by Clifford theory, $\tau\cong \Ind_{K/F}\,\rho$, for some $\rho\in \wW K$. The restriction of $\tau$ to $\scr W_K$ is therefore not irreducible. However, $\scr W_K$ contains $\scr P_K = \scr P_F$ so this contradicts our hypothesis. Part (1) is proved. 
\par 
In part (2), let $\chi\in D(\tau)$. View $\chi$ as a character of $F^\times$ and suppose that $\chi_K = 1$. The field norm $\N KF$ induces a surjection $U^1_K\to U^1_F$, so $\chi$ is trivial on $U^1_F$. Part (1) implies $\chi = 1$. For the second assertion, write $\vG = \Gal KF$. By transitivity and the first assertion, we need only treat the case where $\vG$ is cyclic. Let $\phi\in D(\tau_K)^\vG$. This implies $\phi = \chi_K$, for a character $\chi$ of $F^\times$. We then have $(\chi\otimes\tau)_K \cong \tau_K$, whence it follows that $\chi\otimes\tau \cong \delta\otimes \tau$, for a character $\delta$ of $F^\times$ such that $\delta_K = 1$. In particular, $\delta^{-1}\chi\in D(\tau)$ while $(\delta^{-1}\chi)_K$ agrees with $\phi$ on $U^1_K$. Therefore $(\delta^{-1}\chi)_K = \phi$, as required. \qed 
\enddemo 
\subhead 
4.2 
\endsubhead 
We recall some facts about representations {\it of Heisenberg type.} Let $\tau\in \wW F$ be totally ramified of dimension $p^r$, for some $r\ge1$. Let $\bar\tau:\scr W_F \to \roman{PGL}_{p^r}(\Bbb C)$ be the associated projective representation. Thus $\roman{Ker}\, \bar\tau  = \scr W_E$, for a finite Galois extension $E/F$, and the restriction of $\tau$ to $\scr W_E$ is a multiple of a character $\xi_\tau$ of $\scr W_E$. We call $E$ the {\it p-kernel field\/} of $\tau$ and $\xi_\tau$ the {\it p-central character\/} of $\tau$. 
\par
We are concerned with the case where $\vD = \Gal EF = \roman{Im}\,\bar\tau \cong (\Bbb Z/p\Bbb Z)^{2r}$. The quantity $\xi_\tau[x,y] = \xi_\tau(xyx^{-1}y^{-1})$, $x,y\in \scr W_F$, is then of order dividing $p$. We identify the group $\mu_p(\Bbb C)$ of $p$-th roots of unity in $\Bbb C$ with the additive group of the field $\Bbb F_p$ of $p$ elements. The group $\vD$ is a vector space over $\Bbb F_p$ and the identification $\mu_p(\Bbb C) = \Bbb F_p$ allows us to view the pairing $h_\tau:(x,y)\mapsto \xi_\tau[x,y]$ as a bilinear form $\vD\times \vD \to \Bbb F_p$. This form $h_\tau$ is alternating and nondegenerate. One says that $\tau$ is {\it of Heisenberg type,\/} since the image $\tau(\scr W_F)$ is a Heisenberg group. 
\proclaim{Lemma} 
Let $\tau\in \wW F$ be of Heisenberg type, with p-kernel field $E$, p-central character $\xi_\tau$ and dimension $p^r$.  
\roster 
\item If $\theta$ is an irreducible representation of\/ $\scr W_F$ such that $\theta|_{\scr W_E}$ contains $\xi_\tau$, then $\theta\cong \tau$. 
\item 
Let $\vD'$ be a subgroup of $\vD$, of order $p^r$, such that $h_\tau[x,y]=0$, for all $x,y\in \vD'$. Let $\scr W_{E'}$ be the inverse image of $\vD'$ in $\scr W_F$. 
\itemitem{\rm (a)} There exists a character $\psi$ of\/ $\scr W_{E'}$ such that $\psi|_{\scr W_E} = \xi_\tau$. 
\itemitem{\rm (b)} For any such character $\psi$, the induced representation $\Ind_{E'/F}\,\psi$ is equivalent to $\tau$. In particular, $\Ind_{E/F}\,\xi_\tau$ is a sum of $p^r$ copies of $\tau$. 
\item 
The group $D(\tau)$ has order $p^{2r}$ and consists of all characters of $F^\times$ trivial on norms from $E$. That is, $D(\tau)$ is the group $\widehat\vD$ of characters of $\vD$. 
\endroster 
\endproclaim 
\demo{Proof} 
This is an exercise, for which hints may be found in, for example, \cite{2} 8.3. \qed 
\enddemo 
We recall the structure of the irreducible primitive representations of $\scr W_F$. Leaving aside the trivial one-dimensional case, any such representation has dimens\-ion $p^r$, $r\ge1$, and is totally wildly ramified. 
\proclaim{Proposition } 
Let $\sigma \in \wW F$ be primitive of dimension $p^r$, $r\ge1$. Let $E/F$ be the p-kernel field of $\sigma$, and let $T/F$ be the maximal tamely ramified sub-extension of $E/F$. Let $\vP = \Gal ET$ and $\vG = \Gal TF$. Let $\xi_\sigma$ be the p-central character of $\sigma$.
\roster 
\item  The restriction $\sigma_T$ of $\sigma$ to $\scr W_T$ is irreducible of Heisenberg type. 
\item The alternating form $h_\sigma$ is invariant under $\vG$, 
$$ 
h_\sigma(\gamma u,\gamma v) = h_\sigma(u,v), \quad u,v \in \vP,\ \gamma\in \vG. 
$$ 
\item 
The symplectic $\Bbb F_p$-representation of $\vG$ afforded by $\vP$ is anisotropic and faithful. 
\endroster 
\endproclaim 
\demo{Proof} 
The proposition summarizes \cite{23} Theorem 2.2, Theorem 4.1. \qed 
\enddemo 
We shall refer to the extension $T/F$ as the {\it imprimitivity field\/} of $\sigma$. 
\subhead 
4.3 
\endsubhead 
It will be useful to have an external description of the imprimitivity field $T/F$ arising in 4.2. 
\proclaim{Proposition} 
Let $\sigma \in \wW F$ be primitive of dimension $p^r$. Let $T/F$ be the imprimitivity field of $\sigma$. 
\roster 
\item The group $D(\sigma_T)$ has order $p^{2r}$. 
\item If $L/F$ is a finite tame extension, then $|D(\sigma_L)|\le p^{2r}$. Moreover, $|D(\sigma_L)| = p^{2r}$ if and only if there exists an $F$-embedding $T\to L$. 
\endroster 
In particular, $T/F$ is the unique minimal tame extension for which $D(\sigma_T)$ has order $p^{2r}$. 
\endproclaim 
\demo{Proof} 
Part (1) has been noted in 4.2, and the first assertion of (2) comes from 4.1 Lemma. Next, let $K/F$ be a finite, Galois, tame extension containing both $L$ and $T$. By 4.1 Proposition and part (1), the canonical map $D(\sigma_T) \to D(\sigma_K)$ is bijective while $D(\sigma_L)\to D(\sigma_K)$ is injective. The latter identifies $D(\sigma_L)$ with the group of $\Gal KL$-fixed points in $D(\sigma_K)$. However, if $\phi\in D(\sigma_K)$, the $\Gal KF$-isotropy group of $\phi$ contains $\Gal KT$, whence the result follows. \qed 
\enddemo 
\head \Smc 
5. The Langlands parameter 
\endhead 
Let $\pi\in \upr 1{\Ao{p^r}F}$, for some $r\ge1$, and set $\sigma = \upr L\pi\in \upr 1{\Go{p^r}F}$.  We use the representation $\pi$ to identify the imprimitivity field $T/F$ of $\sigma$, and the  group $D(\sigma_T)$ of characters $\chi$ of $\scr W_T$ such that $\chi\otimes \sigma_T \cong \sigma_T$. As in 4.2 Lemma, the group of characters $D(\sigma_T)$ determines the p-kernel field $E/T$ of $\sigma$. 
\subhead 
5.1 
\endsubhead 
We describe $\pi$, following the outline of 2.1. We therefore choose a character $\psi_F$ of $F$ such that $c(\psi_F) = -1$. To get clean results, it will be necessary to impose the normalization 
$$
\psi_F(\zeta^p) = \psi_F(\zeta),\quad \zeta \in \mu_F, 
\tag 5.1.1 
$$ 
although it will not be invoked until a late stage (5.14 below). 
\par 
There is an m-simple stratum $[\frak a,1,0,\alpha]$ in $A = \M{p^r}F$, with $e_\frak a = p^r$, such that $\pi$ contains the character 
$$
\theta_\alpha: 1{+}x \longmapsto \psi_A(\alpha x), \quad 1{+}x \in U^1_\frak a, 
\tag 5.1.2 
$$
where $\psi_A = \psi_F\circ \roman{tr}_A$. 
\proclaim{Theorem}
The imprimitivity field $T$ of $\sigma$ is the splitting field over $F$ of the polynomial 
$$
X^{p^{2r}-1} - (-1)^p\det\alpha^{p^r-1}. 
\tag 5.1.3 
$$ 
The non-trivial elements of $D(\sigma_T)$ are the characters $\Delta_c$ of $T^\times$ such that 
$$
\Delta_c^p = 1, \quad \Delta_c|_{\mu_T} = 1,\quad \Delta_c(\det\alpha) = 1, 
\tag 5.1.4 
$$ 
and 
$$ 
\Delta_c(1{+}y) = \psi_T(cy),\quad 1{+}y\in U^1_T, 
\tag 5.1.5 
$$
where $c$ ranges over the roots of the polynomial \rom{(5.1.3)} and $\psi_T = \psi_F\circ \roman{Tr}_{T/F}$. 
\endproclaim 
\subhead 
5.2 
\endsubhead 
The proof of the theorem will occupy the rest of the section, but first we draw some conclusions. 
\par 
The conditions (5.1.4), (5.1.5) determine the character $\Delta_c$ uniquely. The representation $\pi$ determines only the coset $(\det\alpha)U^1_F$ (2.2). Changing $\det\alpha$ within that coset changes neither $T$ (as a subfield of $\bar F$) nor the group $D(\sigma_T)$. Let $E/F$ be the p-kernel field of $\sigma$. As in 4.2 Lemma, the abelian extension $E/T$ is given by the relation 
$$ 
\N ET(E^\times) = \bigcap_c \,\roman{Ker}\,\Delta_c, 
\tag 5.2.1 
$$ 
with $c$ ranging over the roots of the polynomial (5.1.3). Thus $E$ is also determined by the coset $(\det\alpha)U^1_F$. 
\par 
In the other direction, the datum $(T/F,D(\sigma_T))$ does {\it not\/} fully determine the coset $(\det\alpha)U^1_F$. 
\proclaim{Corollary} 
Let $\pi_i\in \upr1{\Ao{p^r}F}$ contain the simple stratum $[\frak a,1,0,\alpha_i]$, $i=1,2$. Let $E_i/F$ be the p-kernel field of $\sigma_i = \upr L\pi_i$ and $T_i/F$ its imprimitivity field. The following conditions are equivalent: 
\roster 
\item $E_1 = E_2$; 
\item $T_1 = T_2$; 
\item $\det\alpha_1 \equiv \zeta\det\alpha_2 \pmod{U^1_F}$, for some $\zeta\in F^\times$ such that $\zeta^{p^r} = \zeta$. 
\endroster 
\endproclaim 
\demo{Proof} 
Suppose first that $E_1 = E_2$. Therefore $T_1 = T_2$ and comparison of the polynomials (5.1.3) defining the two fields yields condition (3). Supposing that (3) holds, 5.1 Theorem implies $T_1 = T_2 = T$, say, and $D(\sigma_{1,T}) = D(\sigma_{2,T})$. The relation (5.2.1) implies $E_1 = E_2$, as required for (1). \qed 
\enddemo 
\subhead 
5.3 
\endsubhead 
We start the proof of 5.1 Theorem with a preliminary estimate of the imprimitivity field $T/F$. 
\proclaim{Theorem} 
Let $\sigma\in \upr 1{\Go{p^r}F}$, $r\ge1$, and let $T/F$ be its imprimitivity field. 
\roster 
\item 
The field $T$ satisfies $e(T|F) = 1{+}p^r$ and it contains a root of unity of order $p^{2r}{-}1$. 
\item 
If $\chi\in D(\sigma_T)$ and $\chi\neq 1$, then $\sw(\chi) = 1$. 
\endroster 
\endproclaim 
\demo{Proof} 
The proof takes until the end of the next sub-section. In this one, we find an upper bound for $e(T|F)$. Set $\vP = \Gal ET$ and $\vG = \Gal TF$. 
\proclaim{Proposition} 
The $\Bbb F_p\vG$-representation afforded by $\vP$ is irreducible. 
\endproclaim 
\demo{Proof} 
Suppose otherwise. From \cite{23} Theorem 12.2 we deduce the existence of $\sigma_1,\sigma_2 \in \wW F$ such that $\dim\sigma_i = p^{r_i}$, with $r_2\ge r_1 \ge1$, and $\sigma \cong \sigma_1\otimes \sigma_2$. Let $\pi = \upr L\sigma$, $\pi_i = \upr L\sigma_i$, and let $\ve(\pi_1 \times \pi_2,s,\psi)$ be the local constant of the pair $(\pi_1,\pi_2)$, in the sense of Jacquet, Piatetskii-Shapiro and Shalika \cite{22}, Shahidi \cite{29}. This takes the form 
$$ 
\ve(\pi_1 \times \pi_2,s,\psi) =  q^{-sp^r(\frak F(\check\pi_1,\pi_2)+c(\psi))}\,\ve(\pi_1\times\pi_2,0,\psi), 
$$ 
for a certain rational number $\frak F(\check\pi_1,\pi_2)$ worked out in \cite{6}, \cite{13}. We have 
$$
\ve(\sigma,s,\psi) = \ve(\pi_1\times\pi_2,s,\psi), 
$$ 
so $\sw(\sigma) = p^r\frak F(\check\pi_1,\pi_2) - p^r$. In present circumstances, $\check\pi_1$ cannot be of the form $\chi\pi_2$ for an unramified character $\chi$ of $F^\times$, since that would imply $\sigma_1\otimes\sigma_2$ reducible. So, from \cite{6} 2.1 Proposition and Corollary, we obtain 
$$
\frak F(\check\pi_1,\pi_2) > 1 + \frac{\frak c}{p^{2r_1}}, 
$$
for a certain integer $\frak c \ge 1$. In particular, $\sw(\sigma)>p^{r_2-r_1}\frak c \ge 1$. That is, $\sw(\sigma) > 1$, contrary to hypothesis. \qed 
\enddemo 
\proclaim{Corollary} 
The ramification index $e(T|F)$ divides $1{+}p^r$. 
\endproclaim 
\demo{Proof} 
Let $K/F$ be the maximal unramified sub-extension of $T/F$. The representation $\sigma_K = \sigma|_{\scr W_K}$ satisfies the same hypotheses as $\sigma$: it is irreducible, $\dim\sigma_K = p^r$ and $\sw(\sigma_K) = 1$. It is therefore primitive and the proposition applies to it unchanged. In other words, we may assume that the Galois tame extension $T/F$ is {\it totally ramified.} The cyclic group $\vG = \Gal TF$ thus admits an irreducible, faithful, anisotropic, symplectic $\Bbb F_p$-representation of dimension $2r$. By \cite{10} 3.3, $|\vG|$ divides $1{+}p^r$, as required. \qed 
\enddemo \enddemo 
\remark{Remark} 
The bound $e(T|F)\le 1{+}p^r$ applies to any primitive representation of $\scr W_F$ of dimension $p^r$, irrespective of the value of the Swan conductor \cite{19}. We have used this more demanding technique for its wider interest. 
\endremark 
\subhead 
5.4  
\endsubhead 
We continue the proof of 5.3 Theorem. 
\proclaim{Lemma} 
 Let $K/F$ be a finite tame extension, and set $e_K = e(K|F)$. If there exists a non-trivial character $\chi$ of\/ $\scr W_K$ such that $\chi\otimes\sigma_K \cong \sigma_K$, then $e_K\ge1{+} p^r$. 
\endproclaim 
\demo{Proof} 
The representation $\sigma_K$ satisfies $\sw(\sigma_K) = e_K$. If $\chi$ is a character of $\scr W_K$ with $\sw(\chi) = s\ge0$, then $\sw(\chi\otimes\sigma_K) \le \roman{max}\,\{e_K, p^rs\}$, with equality if $e_K\neq p^rs$ \cite{19} 3.6. A relation $\chi\otimes\sigma_K \cong \sigma_K$ implies first that $s\ge1$ ({\it cf\.} 4.3) and second that $p^rs \le e_K$. Since $p$ does not divide $e_K$, we get $e_K>p^r$, as required. \qed 
\enddemo 
Applying the lemma to the case $K=T$ and recalling 5.3 Corollary, we get $e_T = e(T|F) = 1{+}p^r$. Also, if $\chi\in D(\sigma_T)$ is non-trivial and has Swan conductor $s$, the argument of the preceding proof gives $sp^r<1{+}p^r$, whence $s=1$. 
\par
Finally, since $T/F$ is Galois with ramification index $1{+}p^r$, surely $T$ contains a primitive $1{+}p^r$-th root of unity. That is, if $|\Bbbk_T| = p^t$, say, then $1{+}p^r$ divides $p^t{-}1$. An elementary argument shows that $t$ is divisible by $2r$, whence follows the only remaining assertion of 5.3 Theorem.  \qed 
\subhead 
5.5 
\endsubhead 
We start the proof of 5.1 Theorem. Let $K/F$ be a finite Galois extension with $e(K|F) = 1{+}p^r$: in particular, $K/F$ is tamely ramified. The representation $\sigma_K$ is irreducible, so we define $\pi_K\in \Ao{p^r}K$ by $\upr L\pi_K = \sigma_K$. A character $\chi$ of $K^\times$ then lies in $D(\sigma_K)$ if and only if $\chi\pi_K \cong \pi_K$. We therefore approach 5.1 Theorem via the representation $\pi_K$. Following 4.3 Proposition, we have to find the least extension $K/F$ for which the equation $\chi\pi_K \cong \pi_K$ has $p^{2r}$ solutions $\chi$. 
\par
We describe $\pi_K$ in terms of simple strata and simple characters. The field extension $P = F[\alpha]/F$ is totally wildly ramified, of degree $p^r$. We set $B = \M{p^r}K = K\otimes_FA$. The $K$-algebra $K\otimes_FP$ is a field. We denote it $KP$ and identify it with a subfield of $B$. Let $\frak b$ be the unique hereditary $\frak o_K$-order in $B$ which is stable under conjugation by $KP^\times$. Let $\frak q$ be the Jacobson radical of $\frak b$. The $K$-period $e_\frak b$ of $\frak b$ is $p^r = e(KP|K)$. The quadruple $[\frak b,1{+}p^r,0,\alpha]$ is an m-simple stratum in $B$ and 
$$
H^1(\alpha,\frak b) = U^1_{KP} U_\frak b^{1+[1{+}p^r/2]}, 
$$ 
where $[x]$ denotes the integer part of the real number $x$. 
\par 
We put $\psi_K = \psi_F\circ\roman{Tr}_{K/F}$. Thus $\psi_K$ is a smooth character of $K$ such that $c(\psi_K) = -1$. It satisfies the analogue of (5.1.1), namely 
$$
\psi_K(\zeta^p) = \psi_K(\zeta), \quad \zeta\in \mu_K. 
\tag 5.5.1 
$$ 
Let $\psi_B = \psi_K\circ\roman{tr}_B$. We define a simple character $\theta^K_\alpha\in \scr C(\frak b,\alpha,\psi_K)$ by 
$$
\alignedat3 
\theta^K_\alpha(1{+}h) &= \psi_B(\alpha h),&\quad &1{+}h \in U_\frak b^{1+[1{+}p^r/2]}, \\
\theta^K_\alpha(1{+}y) &= \theta_\alpha(\N{KP}P(1{+}y)), &\quad &1{+}y\in U^1_{KP}. 
\endalignedat 
\tag 5.5.2 
$$
We remark that the character $\theta_\alpha$ is trivial on $U^1_P$, so $\theta^K_\alpha$ is trivial on $U^1_{KP}$. 
\proclaim{Proposition} 
The representation $\pi_K$ contains the m-simple character $\theta^K_\alpha$. 
\endproclaim 
\demo{Proof} 
The character $\theta^K_\alpha$ is m-simple, and its endo-equivalence class is the unique $K/F$-lift of that of $\theta_\alpha$ (see \cite{4} for this concept). The proposition follows from \cite{6} Theorem A. \qed 
\enddemo 
\subhead 
5.6 
\endsubhead 
Let $\phi$ be a character of $K^\times$, $\phi\neq 1$, such that $\phi\pi_K \cong \pi_K$. Viewing $\phi$ as a character of $\scr W_K$ via class field theory, we get $\phi\otimes \sigma_K \cong \sigma_K$. We argue as in 5.4 to conclude that $\roman{sw}(\phi) = 1$. Consequently,  there is an element $c\in K$, with $\ups_K(c) = -1$, such that $\phi(1{+}y) = \psi_K(cy)$, $y\in \frak p_K$. 
\proclaim{Proposition} 
Let $c\in K$, $\ups_K(c) = -1$. 
\roster 
\item The quadruple $[\frak b,1{+}p^r,0,\alpha{+}c]$ is an m-simple stratum in $B$ and 
$$ 
H^1(\alpha{+}c,\frak b) = H^1(\alpha,\frak b). 
$$
\item 
Let $\phi$ be a character of $K^\times$ such that $\phi(1{+}y) = \psi_K(cy)$, $y\in \frak p_K$. Let $\vF$ denote the character $\phi\circ\det_B|_{H^1(\alpha,\frak b)}$. 
\itemitem{\rm (a)} The character $\vF\theta^K_\alpha$ lies in $\scr C(\frak b,\alpha{+}c,\psi_K)$ and is contained in $\phi\pi_K$. 
\itemitem{\rm (b)} If $\phi\pi_K \cong \pi_K$, the simple characters $\theta_\alpha^K$, $\vF\theta_\alpha^K$ are conjugate in $\GL{p^r}K$. 
\itemitem{\rm (c)} If the simple characters $\theta_\alpha^K$, $\vF\theta_\alpha^K$ are conjugate in $\GL{p^r}K$, there is a unique character $\phi'$ of $K^\times$ that agrees with $\phi$ on $U^1_K$ and such that $\phi'\pi_K \cong \pi_K$. 
\endroster 
\endproclaim 
\demo{Proof} 
Assertion (1) is immediate from the definitions. Part (2)(a) is an instance of \cite{15} Appendix. Part (2)(b) is given by \cite{14} 8.4 (or \cite{11} Corollary 1 for this formulation). In part (2)(c), the totally ramified representations $\pi_K$, $\phi\pi_K$ contain the same m-simple character and so $\phi\pi_K \cong \chi\pi_K$, for a tamely ramified character $\chi$ of $K^\times$. The result therefore holds with $\phi' = \chi^{-1}\phi$. The uniqueness of $\phi'$ is given by 4.1 Proposition. \qed 
\enddemo  
\remark{Remark} 
In the context of part (2)(a) of the proposition, $\scr C(\frak b,\alpha{+}c,\psi_K) = \vF\,\scr C(\frak b,\alpha,\psi_K)$. 
\endremark 
\subhead 
5.7 
\endsubhead 
We collect some technical results. 
\proclaim{Lemma 1} 
If $d(KP|K)$ is the differental exponent of $KP/K$, then $d(KP|K) \ge 2p^r$. 
\endproclaim 
\demo{Proof} 
The extension $P/F$ is totally wildly ramified, so $d(P|F) \ge [P{:}F] = p^r$. The transitivity property of the different gives us first 
$$
d(KP|F) = d(KP|P)+(1{+}p^r)d(P|F) \ge p^r+(1{+}p^r)p^r = 2p^r+p^{2r}, 
$$ 
and second $d(KP|F) = d(KP|K)+p^{2r}$. It follows that $d(KP|K) \ge 2p^r$, as required. \qed 
\enddemo 
Let $s_K:B\to KP$ be a tame corestriction on $B$, relative to $KP/K$ (see \cite{14} (1.3.3) for the definition). 
\proclaim{Lemma 2} 
For any integer $t$, we have $s_K(\frak p_{KP}^t) \subset \frak p_{KP}^{1+p^r+t}$. 
\endproclaim 
\demo{Proof} 
This follows from Lemma 1 and \cite{14} (1.3.8). \qed 
\enddemo 
\subhead 
5.8 
\endsubhead 
We take $c\in K$ with $\ups_K(c) = -1$. We compare the m-simple strata $[\frak b,1{+}p^r,0,\alpha]$, $[\frak b,1{+}p^r,0,\alpha{+}c]$. 
\proclaim{Proposition} 
There exists $x\in \frak q$ such that 
$$
(1{+}x)^{-1}\alpha (1{+}x) \equiv \alpha{+}c \pmod{\frak b}. 
\tag 5.8.1 
$$
\endproclaim 
\demo{Proof} 
Let $\Bbb A:B \to B$ be the map $x \mapsto \alpha x\alpha^{-1}-x$. This fits into a long exact sequence \cite{14} 1.4  
$$
\cdots \to B @>{\ \Bbb A\ }>> B @>{\ s_K\ }>> B@>{\ \Bbb A\ }>> B \to \cdots . 
$$ 
Set $\delta  = c^{-1}\alpha$: thus $\delta^{-1}$ is a prime element of $KP$ and $\delta^{-1}\frak b = \frak q$. Multiplying the congruence (5.8.1) on the left by $c^{-1}(1{+}x)$, we see it is equivalent to 
$$ 
\Bbb A(x)-x\delta^{-1} \equiv \delta^{-1} \pmod{\frak q^{p^r}}. 
\tag 5.8.2 
$$ 
\indent 
Let $i,j$ be integers, with $1\le j\le p^r$. Since $\frak p_K\frak b = \frak q^{p^r}$, the quotient $V_{i,j} = \frak q^i/\frak q^{i+j}$ is a $\Bbbk_K$-vector space. Each of the maps $s_K$, $\Bbb A$ induces a $\Bbbk_K$-endomorphism of $V_{i,j}$, for which we use the same notation. Since $\alpha$ is minimal over $K$, we have an infinite exact sequence \cite{14} (1.4.7), (1.4.15), 
$$
\cdots \to V_{i,j} @>{\Bbb A}>> V_{i,j} @>{s_K}>> V_{i,j} @>{\Bbb A}>> V_{i,j} \to \cdots . 
\tag 5.8.3 
$$ 
The subspace $s_K(V_{i,j})$ is the natural image of $\frak p_{KP}^i/\frak p_{KP}^{i+j}$ in $V_{i,j}$. 
\proclaim{Lemma} 
\roster 
\item  
The endomorphism $\Bbb A$ of $V_{i,j}$ is nilpotent such that $\Bbb A^{p^r} = 0$ and 
$$ 
\Bbb A^{p^r-1}(V_{i,j}) = s_K(V_{i,j}) = \frak p_{KP}^i/\frak p_{KP}^{i+j} \neq 0. 
$$ 
\item There is a unit $u\in U_{KP}$ such that 
$$ 
\Bbb A^{p^r-1}(v) = us_K(v), \quad v\in V_{i,j}. 
$$ 
\endroster 
\endproclaim 
\demo{Proof} 
We consider the relation between the hereditary orders $\frak a$, $\frak b$. Let $\scr L$ be the chain of $\frak a$-lattices in the vector space $F^{(p^r)}$. There is then a unique $\frak o_K$-lattice chain in $K\otimes F^{(p^r)}$ which is stable under translation by elements of $(K\otimes P)^\times = KP^\times$, and contains every lattice $\frak o_K\otimes_{\frak o_F} L$, $L\in \scr L$. The hereditary $\frak o_K$-order defined by this chain is $\frak b$. Consequently, if $\frak p = \roman{rad}\,\frak a$, then $\frak p \subset \frak q^{1+p^r}$. 
\par 
The element $\alpha^{p^r}$ lies in $\vp_\alpha^{-1}U^1_\frak a\subset \vp^{-1}_\alpha U_\frak b^{1+p^r}$, for a certain prime element $\vp_\alpha$ of $F$. Thus $\alpha^{p^r}$ acts trivially on $V_{i,j}$, and so $\Bbb A^{p^r}(V_{i,j}) = 0$. The exact sequence (5.8.3) shows that $\roman{Ker}\,\Bbb A$ has dimension $j = p^{-r}\dim V_{i,j}$. We conclude that $\Bbb A$ is nilpotent of type $(j,j,\dots,j)$ and all assertions in (1) follow. 
\par 
It follows that $\roman{Ker}\,\Bbb A^{p^r-1} = \roman{Im}\,\Bbb A = \roman{Ker}\,s_K$. The maps $s_K$ and $\Bbb A^{p^r-1}$ are therefore $\frak o_{KP}$-isomorphisms between two free $\frak o_{KP}/\frak p_{KP}^j$-modules of rank one, whence they differ by a unit. \qed 
\enddemo 
We apply the lemma to the case $V = V_{1,p^r-1} = \frak q/\frak q^{p^r}$. We define another $\Bbbk_K$-endomorphism $\Bbb B_c$ of $V$ by 
$$
\Bbb B_c(v) = \Bbb A(v)-v\delta^{-1}, \quad v\in V. 
$$
To prove the proposition, we have to show that the equation $\Bbb B_c(x) = \delta^{-1}$ has a solution $x = x_c$ in $V$. To do this, we abbreviate $m=p^r{-}2$ and define 
$$ 
\multline 
X_c(z) = \Bbb A^m(z)+\Bbb A^{m-1}(z)\delta^{-1} + \Bbb A^{m-2}(z)\delta^{-2} + \dots \\ 
\dots + \Bbb A(z)\delta^{1-m} +z\delta^{-m}, \quad z\in V, 
\endmultline 
$$ 
so that 
$$
\Bbb B_c(X_c(z)) = \Bbb A^{m+1}(z)-z\delta^{-(m+1)} = \Bbb A^{m+1}(z) 
$$
in $V$. By the lemma, we may choose $z\in V$ such that $\Bbb A^{m+1}(z) = \Bbb A^{p^r-1}(z) = \delta^{-1}$, and then $x = X_c(z)$ provides the desired solution to the congruence (5.8.2). \qed 
\enddemo 
\subhead 
5.9 
\endsubhead 
We gather the threads for the next phase of the proof of 5.1 Theorem. 
\par 
Let $c\in K$, $\ups_K(c) = -1$ and let $x = x_c\in \frak q$ satisfy 
$$ 
(1{+}x)^{-1}\alpha (1{+}x) \equiv \alpha{+}c \pmod{\frak b},  
$$
as in 5.8 Proposition. The formula $\psi_{K,c}:1{+}t\mapsto \psi_K(ct)$ defines a character of $U^1_K$. This gives a character 
$$ 
\chi_c:h\longmapsto \psi_{K,c}({\det}_B h) 
$$ 
of $H^1(\alpha,\frak b)$. As in 5.6 Remark, $\scr C(\frak b,\alpha{+}c,\psi_K) = \scr C(\frak b,\alpha,\psi_K)\,\chi_c$. 
\proclaim{Lemma} 
The element $1{+}x$ normalizes $H^1(\alpha,\frak b) = H^1(\alpha{+}c,\frak b)$ and conjugation by $1{+}x$ gives a bijection 
$$
\align 
\scr C(\frak b,\alpha,\psi_K) &\longrightarrow \scr C(\frak b,\alpha{+}c,\psi_K), \\ \vt &\longmapsto \vt^{1+x}. 
\endalign 
$$ 
There exists $\vt_c\in \scr C(\frak b,\alpha,\psi_K)$ such that $(\theta^K_\alpha)^{1+x} = \vt_c\chi_c$. 
\endproclaim 
\demo{Proof} 
The equality of $H^1$-groups has been noted in 5.6 Proposition, whence the other assertions follow. \qed 
\enddemo 
The next step in the proof of 5.1 Theorem is to find all elements $c$ for which $\vt_c = \theta^K_\alpha$, that is, 
$$
\big(\theta^K_\alpha\big)^{1+x_c} = \theta^K_\alpha\chi_c. 
\tag 5.9.1 
$$ 
For such an element $c$, 5.6 Proposition gives a character $\Delta_c$ of $K^\times$, extending $\psi_{K,c}$, such that $\Delta_c \pi_K \cong \pi_K$. 
\subhead 
5.10 
\endsubhead 
We have to compute the quantity 
$$ 
\theta^K_\alpha((1{+}x)h(1{+}x)^{-1}), \quad h\in H^1(\alpha,\frak b) = U^1_{KP}U^{1+[1+p^r/2]}_\frak b, 
$$ 
where $x = x_c$. The definition of $x_c$ yields 
$$
\theta^K_\alpha((1{+}x)h(1{+}x)^{-1}) = \theta^K_\alpha(h)\chi_c(h), \quad h \in U^{1+[1+p^r/2]}_\frak b. 
\tag 5.10.1 
$$
It is therefore enough to consider elements $h = 1{+}y\in U^1_{KP}$. 
\proclaim{Proposition} 
If $y\in \frak p_{KP}$ and $x = x_c$, then $(\theta^K_\alpha)^{1+x}(1{+}y) = \psi_B(-cxy)$. 
\endproclaim 
\demo{Proof} 
We need a preliminary calculation. 
\proclaim{Lemma} 
If $y\in \frak p_{KP}$, then 
$$
(1{+}x)(1{+}y)(1{+}x)^{-1} \equiv 1{+}v \pmod{U^{1+p^r}_\frak b}, 
$$
where $v\in \frak p_{KP}$ and $v\equiv y \pmod{U^2_{KP}}$. 
\endproclaim 
\demo{Proof}
We write $[x,y] = xy{-}yx$, so that 
$$
1{+}t = (1{+}x)(1{+}y)(1{+}x)^{-1} = 1+y+[x,y](1{+}x)^{-1}. 
$$ 
Surely $t\equiv y\pmod{\frak q^2}$, so the result reduces to showing $\alpha t{-}t\alpha\in \alpha\frak q^{1+p^r} = \frak b$. We use two forms of the defining relation for $x$, namely 
$$
\align 
\alpha x{-}x\alpha &\equiv c(1{+}x) \pmod{\frak b}, \\ 
(1{+}x)^{-1}\alpha &\equiv (\alpha{+}c)(1{+}x)^{-1} \pmod{\frak b}. 
\endalign 
$$ 
Expanding, 
$$ 
\align
\alpha t{-}t\alpha &= \alpha[x,y](1{+}x)^{-1}- [x,y](1{+}x)^{-1}\alpha \\ 
&= \big(\alpha [x,y] - [x,y](\alpha{+}c)\big)(1{+}x)^{-1} \\ 
&\equiv \alpha [x,y] - [x,y](\alpha{+}c) \pmod{\frak b}. 
\endalign 
$$ 
We recall that $c$ is central and that $y$ commutes with $\alpha$. Expanding the commutators, the defining relation yields $
\alpha t{-}t\alpha \equiv 0 \pmod{\frak b}$, as required. \qed 
\enddemo  
We use the lemma to write 
$$ 
(1{+}x)(1{+}y)(1{+}x)^{-1} = 1{+}v{+}h, \quad h\in \frak q^{1+p^r},\ v\in \frak p_{KP}, 
$$
with $v\equiv y\pmod{\frak p_{KP}^2}$. Thus 
$$ 
(\theta^K_\alpha)^{1+x}(1{+}y) = \theta^K_\alpha(1{+}y)\,\psi_B(\alpha h) .  
$$ 
We have $\theta^K_\alpha(1{+}y) = 1 = \psi_B(\alpha y)$ by (5.5.2), 5.7 Lemma 2, respectively, so 
$$
(\theta^K_\alpha)^{1+x}(1{+}y) = \psi_B(\alpha(v{+}h)), \quad\text{and}\quad v{+}h = y+[x,y](1{+}x)^{-1}. 
$$ 
We have to compute 
$$
\psi_B(\alpha [x,y](1{+}x)^{-1}) = \psi_B\big(\sum_{i\ge0} (-1)^i \alpha(xy{-}yx)x^i\big). 
$$ 
We expand the inner sum using the symmetry properties of the trace and the defining relation for $x$ in the form 
$$
[\alpha, x] \equiv c(1{+}x) \pmod{\frak b}. 
$$ 
The term $i=0$ contributes $\psi_B(\alpha(xy{-}yx)) = 1$, since $y$ commutes with $\alpha$. The general term $i\ge1$ gives 
$$
\psi_B\big((-1)^i \alpha(xy{-}yx)x^i\big) = \psi_B\big((-1)^ix^i[\alpha,x]y\big) = \psi_B\big((-1)^i c(x^i{+}x^{i+1})\big). 
$$
In all, we get 
$$\psi_B\big(\alpha [x,y](1{+}x)^{-1}\big) = \psi_B(-cxy), 
$$ 
or $(\theta^K_\alpha)^{1+x}(1{+}y) = \psi_B(-cxy)$, as required. \qed 
\enddemo 
\subhead 
5.11 
\endsubhead 
We examine the character 
$$
\eta_c:1{+}y\longmapsto \psi_B(-cx_cy) 
$$
of $U^1_{KP}$ arising in 5.10 Proposition. We use the map $\Bbb A$ of 5.8 and let $s^K_\alpha$ be a tame corestriction on $B$, relative to $KP/K$, such that 
$$
s^K_\alpha(t) \equiv \Bbb A^{p^r-1}(t) \pmod{\frak q^{p^r}}, \quad t\in \frak q, 
$$
({\it cf\.} 5.8 Lemma (2)). 
\proclaim{Lemma} 
There is a unique character $\eps$ of $KP$ such that $c(\eps) = -1$ and 
$$ 
\psi_B(b) = \eps(s^K_\alpha(b)), \quad b\in B. 
\tag 5.11.1
$$ 
The character $\eta_c$ then satisfies 
$$
\eta_c(1{+}y) = \eps(-\alpha^{1-p^r}c^{p^r}y), \quad y\in \frak p_{KP}. 
\tag 5.11.2 
$$
Consequently, $\eta_c$ is trivial on $U^2_{KP}$. 
\endproclaim 
\demo{Proof} 
The first assertion is given by \cite{14} (1.3.7). 
\par 
We return to the construction of $x_c$ in 5.8. As there, we put $\delta = c^{-1}\alpha$. We choose $z\in \frak q$ such that $\Bbb A^{p^r-1}(z) \equiv \delta^{-1} \pmod{\frak q^{p^r}}$: this yields $x_c = X_c(z)$. The choice of $s^K_\alpha$ then gives 
$$
s^K_\alpha(cx_c) \equiv c\delta^{1-p^r} \pmod{\frak b}, 
$$
while $c\delta^{1-p^r} = \alpha^{1-p^r}c^{p^r}$. The lemma now follows from the relation (5.11.1). \qed 
\enddemo 
\subhead 
5.12 
\endsubhead 
We compare the characters $\eta_c$, $\chi_c$ on $U^1_{KP}$. Both are trivial on $U^2_{KP}$ and, for $y\in \frak p_{KP}$, 
$$
\chi_c(1{+}y) = \psi_{K,c}(\N{KP}K(1{+}y)). 
$$
We may take $y = \zeta c\alpha^{-1}$, for some $\zeta\in \mu_{KP} = \mu_K$. 
\proclaim{Lemma} 
If $y = \zeta c\alpha^{-1}$, $\zeta\in \mu_K$, then 
$$ 
\N{KP}K(1{+}y) \equiv 1+\zeta^{p^r}c^{p^r}\det\alpha^{-1} \pmod{U^2_K}. 
$$
\endproclaim 
\demo{Proof} 
Nothing is changed if we replace $\alpha$ by $\alpha u$, $u\in U^1_\frak a$. We may therefore assume (as in the proof of 2.2 Proposition) that $\alpha^{-p^r}$ is a prime element of $F$. Therefore $(c\alpha^{-1})^{p^r}$ is a prime element of $K$. We may choose a matrix representation so that $c\alpha^{-1}$ is a monomial matrix as in 2.2. We may re-write the identity to be proved as 
$$
\det(1{+}y) \equiv 1+ \det(\zeta c\alpha^{-1}) \pmod{U^2_K}. 
$$ 
The lemma is then given by an elementary calculation. \qed 
\enddemo 
It follows that 
$$
\chi_c(1{+}y) = \psi_K(\zeta^{p^r}c^{1+p^r}\det\alpha^{-1}), \quad y = \zeta c\alpha^{-1}. 
\tag 5.12.1 
$$
In these terms, (5.11.2) says 
$$ 
\eta_c(1{+}y) = \eps((-1)^{p^r}c^{1+p^r}\zeta\det\alpha^{-1}), \quad y = \zeta c\alpha^{-1}. 
\tag 5.12.2 
$$
\subhead 
5.13 
\endsubhead 
We relate the characters $\psi_K$, $\eps$. 
\proclaim{Lemma} 
The characters $\eps$, $\psi_K$ satisfy $\eps(\zeta) = \psi_K(\zeta)$, for all $\zeta\in \mu_K$. 
\endproclaim 
\demo{Proof} 
We continue in the situation of the proof of 5.12 Lemma. For $\zeta\in \mu_K$, we let $\zeta_0$ be the matrix with $\zeta$ in the $(1,1)$-place and with all other entries zero. The matrix $\Bbb A^{p^r-1}(\zeta_0)$ is a scalar matrix whose last entry is $\zeta$, whence the lemma follows. \qed 
\enddemo 
\subhead 
5.14 
\endsubhead 
Let $\vs_c\in \mu_K$ satisfy $\vs_c \equiv c^{1+p^r}\det\alpha^{-1} \pmod{U^1_{KP}}$. Combining (5.12.1), (5.12.2), we have to solve the  equation 
$$
\psi_K((-1)^{1+p^r}\vs_c\zeta) = \psi_K(\vs_c\zeta^{p^r}), \quad \zeta \in \mu_K,  
$$ 
for $c\in K$, $\ups_K(c)  = -1$. We have $\psi_K(\gamma^p) = \psi_K(\gamma)$, for $\gamma\in \mu_K$ (5.5.1). We must therefore solve $(-1)^{p^{2r}}\vs_c^{p^r} = \vs_c$, that is, $\vs_c^{1-p^r} = (-1)^p$. This translates back to 
$$
c^{1-p^{2r}} \equiv (-1)^p\det\alpha^{1-p^r} \pmod{U^1_K}. 
$$ 
This congruence admits $p^{2r}{-}1$ solutions in $K^\times/U^1_K$ if and only if the field $K$ contains the splitting field of the polynomial (5.1.3). 
\par
Let $K$ be this splitting field and let $c\in K$ be a root of the polynomial (5.1.3). Let $\phi_c$ be a character of $K^\times$ such that $\phi_c(1{+}t) = \psi_K(ct)$, $t\in \frak p_K$. The representations $\pi_K$, $\phi_c\pi_K$ both contain the simple character $\theta^K_\alpha$. As in 5.6 Proposition, there is a unique character $\Delta_c$ of $K^\times$, agreeing with $\phi_c$ on $U^1_K$, such that $\Delta_c\pi_K \cong\pi_K$. That is, $\Delta_c\in D(\pi_K)$. Surely $K$ is minimal for the property $|D(\sigma_K)| = p^{2r}$, so $K = T$, as desired. 
\par
This character $\Delta_c$ satisfies (5.1.5) by construction. By (5.1.5), $\Delta_c^p$ is tamely ramified and so trivial, by 4.1 Proposition. The second property in (5.1.4) expresses the fact that $\pi_K$, $\Delta_c\pi_K$ have the same central character, while the third follows from the relation $\ve(\Delta_c\pi_K,s,\psi_T) = \ve(\pi_K,s,\psi_T)$. 
\par 
We have completed the proof of 5.1 Theorem. \qed 
\head\Smc 
6. The p-central character 
\endhead 
As in \S5, $\pi \in \upr 1{\Ao {p^r}F}$ contains the simple stratum $[\frak a,1,0,\alpha]$, relative to the character $\psi_F$ of (5.1.1). From 5.1 Theorem, the primitive representation $\sigma = \upr L\pi$ has p-kernel field $E/F$ and imprimitivity field $T/F$. We investigate the p-central character $\xi_\sigma$ of $\sigma$ to complete our examination of the representation $\sigma$. 
\subhead 
6.1 
\endsubhead 
Let $F'/F$ be the maximal unramified sub-extension of $T/F$. Thus $T/F'$ is cyclic and totally ramified of degree $1{+}p^r$. We 
set 
$$ 
\aligned 
\vP &= \Gal ET \cong (\Bbb Z/p\Bbb Z)^{2r}, \\ \vG &= \Gal T{F'} \cong \Bbb Z/(1{+}p^r)\Bbb Z, \\ \vT &= \Gal E{F'}. 
\endaligned \tag 6.1.1 
$$ 
The extension $E/T$ is totally wildly ramified, so there exists $\beta_E\in E^\times$ such that $\N ET(\beta_E) \equiv \det\alpha \pmod{U^1_T}$. This condition determines the coset $\beta_E U^1_E$ uniquely. We write $\psi_E = \psi_F\circ \roman{Tr}_{E/F}$. 
\proclaim{Theorem} 
The character $\xi_\sigma$ has the following properties: 
\roster 
\item $\xi_\sigma$ is fixed under conjugation by $\vT$, 
\item $\sw(\xi_\sigma) = 1{+}p^r$, and 
\item 
if $\beta_E\in E$ satisfies $\N ET(\beta_E) \equiv \det\alpha \pmod{U^1_T}$, then 
$$ 
\xi_\sigma(1{+}x) = \psi_E(\beta_E^{p^r} x), \quad 1{+}x\in U^{1+p^r}_E. 
$$ 
\endroster 
The properties \rom{(1)--(3)} determine the character $\xi_\sigma|_{U^1_E}$ uniquely. 
\endproclaim 
We give a more complete formula for $\xi_\sigma|_{U^1_E}$ in (6.7.2) below. We note the following refinement of 5.2 Corollary. 
\proclaim{Corollary} 
Let $\pi_i \in \upr1{\Ao{p^r}F}$ contain the m-simple stratum $[\frak a,1,0, \allowmathbreak \alpha_i]$, $i=1,2$. Suppose that the representations $\sigma_i = \upr L
\pi_i$ have p-kernel field $E$. If $\xi_i$ is the p-central character of $\sigma_i$, the following conditions are 
equivalent. 
\roster 
\item $\xi_1|_{U^1_E} = \xi_2|_{U^1_E}$; 
\item $\xi_1|_{U^{1+p^r}_E} = \xi_2|_{U^{1+p^r}_E}$; 
\item $\det\alpha_1 \equiv \det\alpha_2 \pmod{U^1_F}$. 
\endroster 
\endproclaim  
The proofs occupy the rest of the section. 
\subhead 
6.2 
\endsubhead 
Before starting, we show how 6.1 Theorem completes our picture of the representation $\sigma$. We use the viewpoint of \cite{12} \S1. According to the Ramification Theorem of local class field theory, the Artin Reciprocity map $\scr W_E \to E^\times$ maps the wild inertia subgroup $\scr P_E$ onto $U^1_E$. From this point of view, the restricted character $\xi_\sigma|_{U^1_E}$, determined in 6.1 Theorem, gives a character $\zeta_\sigma$ of $\scr P_E$. (If we think of $\xi_\sigma$ as a character of $\scr W_E$, then $\zeta_\sigma = \xi_\sigma|_{\scr P_E}$.) 
\proclaim{Corollary} 
\roster 
\item The representation $\rho_\sigma = \sigma|_{\scr P_F}$ is irreducible. It is the unique irreducible representation of $\scr P_F$ containing $\zeta_\sigma$. 
\item 
The representation $\sigma$ is the unique element of $\upr1{\Go{p^r}F}$ with the following properties: 
\itemitem{\rm (a)} $\sigma|_{\scr P_F} \cong \rho_\sigma$, 
\itemitem{\rm (b)} $\det\sigma = \omega_\pi$, and 
\itemitem{\rm (c)} $\ve(\sigma,\frac12,\psi_F) = \ve(\pi,\frac12,\psi_F)$. 
\endroster 
\endproclaim 
\demo{Proof} 
The representation $\rho_\sigma$ is surely irreducible and contains $\zeta_\sigma$. The relation $\Ind_{\scr W_E}^{\scr W_T}\,\xi_\sigma = p^r\sigma_T$ of 4.2 Lemma implies $\Ind_{\scr P_E}^{\scr P_F}\, \zeta_\sigma = p^r\rho_\sigma$, whence (1) follows. 
\par 
That $\sigma$ has the listed properties follows from the definition of $\rho_\sigma$ and 2.3 Proposition. If $\sigma'$ is some other representation of $\scr W_F$ extending $\rho_\sigma$, 1.3 Proposition of \cite{12} asserts that $\sigma'\cong \chi\otimes\sigma$, for a unique tamely ramified character $\chi$ of $\scr W_F$. If $\det\sigma' = \chi^{p^r}\det\sigma = \det\sigma$, then $\chi^{p^r} = 1$ and consequently $\chi$ is unramified. So, by 2.3 Lemma and Proposition, 
$$
\ve(\sigma',\tfrac12,\psi_F) = \chi(\det\alpha)^{-1}\,\ve(\sigma,\tfrac12,\psi_F). 
$$
Since $\ups_F(\det\alpha) = -1$, a relation $\ve(\sigma',\tfrac12,\psi_F) = \ve(\sigma,\tfrac12,\psi_F)$ implies $\chi=1$. \qed 
\enddemo 
\remark{Remark} 
The field extension $E/F$ and the representation $\rho_\sigma$ of $\scr P_F$ are determined entirely by the coset $(\det\alpha)U^1_F$ or, equivalently, by the endo-class $\vT_\alpha$ of the simple character $\theta_\alpha$ occurring in $\pi$ ({\it cf\.} 2.1). Indeed, $\vT_\alpha$ is the endo-class corresponding to $\rho_\sigma$ via the Ramification Theorem of \cite{12} 6.1. 
\endremark 
\subhead 
6.3 
\endsubhead 
We use the notation set up at start of 6.1. 
\par 
A subgroup $\vX$ of $\vP$ is {\it $\xi_\sigma$-Lagrangian\/} if $|\vX| = p^r$ and the commutator pairing $(x,y) \mapsto \xi_\sigma[x,y]$ is null on $\vX$. Equivalently, $\vX$ is the image in $\vP$ of a maximal abelian subgroup $\widetilde\vX$ of $\sigma(\scr W_T)$. It follows that $\xi_\sigma$ extends to a character of $\widetilde\vX$ or, in terms of fields, $\xi_\sigma$ factors through the norm map $\N EK:E^\times \to K^\times$ where $K = E^\vX$, {\it cf\.} 4.2 Lemma. We identify the dual $\widehat\vX$ of the abelian group $\vX$ with the group of characters of $K^\times$ vanishing on norms from $E$. 
\proclaim{Proposition} 
Let $\vX$ be an $\xi_\sigma$-Lagrangian subgroup of $\vP$ and set $K = E^\vX$. 
\roster 
\item 
The fields $T\subset K\subset E$ satisfy 
$$ \gather 
d(E|T) = 2p^{2r}{-}2, \\
d(E|K) = d(K|T) = 2p^r{-}2. 
\endgather 
$$
\item 
If $\phi \in \widehat\vX$, $\phi\neq 1$, then $\sw(\phi) = 1$. 
\item 
The character $\xi_\sigma$ satisfies $\sw(\xi_\sigma) = 1{+}p^r$. 
\item 
If $\chi$ is a character of $K^\times$, such that $\xi_\sigma = \chi\circ \N EK$, then $\sw(\chi) = 2$. 
\endroster 
\endproclaim 
\demo{Proof} 
Each of the extensions $E/K$, $K/T$ is totally ramified and $\vX \cong (\Bbb Z/p\Bbb Z)^r \cong \Gal KT$. As in 5.1, every non-trivial $\phi\in D(\sigma_T) = \widehat\vP$ satisfies $\sw(\phi) = 1$, so the Artin conductor $\roman a(\phi)$ is $2$. By the conductor-discriminant formula \cite{28} VI \S3 Corollaire 2, 
$$
d(E|T) = \sum_{\phi\in D(\sigma_T)} \roman a(\phi) = 2p^{2r}{-}2. 
$$ 
Let $\vPs = \Gal KT$. The dual $\widehat\vPs$ of $\vPs$ is a subgroup of $D(\sigma_T) = \widehat\vP$ of order $p^r$. The relation $d(K|T) = 2p^r{-}2$ follows as before. Multiplicativity of the different in towers yields the final part of (1). Part (2) now follows from the conductor-discriminant formula. 
\par
In part (3), we recall that $\Ind_{E/T}\,\xi_\sigma = p^r\sigma_T$. Since $\sw(\sigma_T) = 1{+}p^r$, the result is given by part (1) and (1.3.1). Finally, $\Ind_{K/T}\,\chi = \sigma_T$ (4.2 Lemma) and part (4) follows similarly. \qed 
\enddemo  
\subhead 
6.4 
\endsubhead 
We take $\vX$ and $K = E^\vX$ as before. 
\proclaim{Proposition} 
The norm map $\N EK$ induces an isomorphism 
$$ 
\align 
U^{1+p^r}_E/U^{2+p^r}_E &\cong U^2_K/U^3_K , \\ \N EK(1{+}x) &\equiv 1+\roman{Tr}_{E/K}(x) \pmod{U^3_K}, 
\endalign 
$$ 
and an exact sequence 
$$ 
1\to \vX \longrightarrow U^1_E/U^2_E @>{\ \N EK}>> U^1_K/U^2_K. 
$$
\endproclaim 
\demo{Proof} 
We consider the ramification subgroups $\vX_i$, $i\in \Bbb Z$, $i\ge0$, of $\vX = \Gal EK$ in the lower numbering. We recall from \cite{28} IV \S2 Proposition 4 that 
$$ 
2p^r{-}2 = d(E|K) = \sum_{i\ge0} \big(|\vX_i|-1\big). 
$$
Here we have $\vX  = \vX_0 = \vX_1$ and $|\vX| = p^r$. The only conclusion is that $\vX_2$ is trivial. The derivative of the Herbrand function $\vf_{E/K}$ is therefore given by 
$$
\varphi'_{E/K}(x) = \left\{\,\alignedat3 &1, &\quad &0<x<1, \\ &p^{-r}, &\quad &1<x, \endalignedat \right. 
$$ 
so the inverse Herbrand function $\varPsi_{E/K}$ satisfies 
$$
\varPsi_{E/K}'(x) =  \left\{\,\alignedat3 &1, &\quad &0<x<1, \\ &p^r, &\quad &1<x. \endalignedat \right. 
$$ 
Now we use \cite{28} V Proposition 9. Abbreviating $\varPsi = \varPsi_{E/K}$, we have an exact sequence 
$$
1 \to \vX_{\varPsi(i)}/\vX_{\varPsi(i)+1} \longrightarrow U_E^{\varPsi(i)}/U_E^{\varPsi(i)+1} @>{\ \N EK\ }>> U^i_K/U^{i+1}_K, 
\tag 6.4.1 
$$
for any integer $i\ge1$. 
\par
Since $\varPsi(2) = 1{+}p^r$, the norm $\N EK$ induces an isomorphism $U^{1+p^r}_E/U^{2+p^r}_E \cong U^2_K/U^3_K$. Since $d(E|K) = 2p^r{-}2$ (6.3), the trace $\roman{Tr}_{E/K}$ induces an isomorphism $\frak p_E^{1+p^r}/\frak p_E^{2+p^r} \to \frak p_K^2/\frak p_K^3$. On the other hand, for $x\in \frak p_E^{1+p^r}$, 
$$ 
\N EK(1{+}x) = 1+\roman{Tr}_{E/K}(x) + R, 
$$ 
where the remainder term $R$ lies in $\frak p_E^{2+2p^r}\cap K \subset \frak p_K^3$. 
\par 
In the second assertion, the exact sequence is the case $i=1$ of (6.4.1). \qed 
\enddemo 
\subhead 
6.5 
\endsubhead 
Again let $\vX$ be an $\xi_\sigma$-Lagrangian subgroup of $\vP$ and $K = E^\vX$. 
\proclaim{Proposition} 
Let $\chi$ be a character of $K^\times$ such that $\xi_\sigma = \chi\circ\N EK$. Let $\psi_K = \psi_F\circ \roman{Tr}_{K/F}$ and let $\beta_K \in K^\times$. The following conditions are equivalent: 
\roster 
\item 
$\chi(1{+}y) = \psi_K(\beta_Ky)$, $y\in \frak p_K^2$; 
\item  
$\N KT\,\beta_K \equiv \det \alpha \pmod{U^1_T}$. 
\endroster 
\endproclaim 
\demo{Proof} 
If $\psi_T = \psi_F\circ\roman{Tr}_{T/F}$, then $c(\psi_T) = -1$ since $T/F$ is tame. It follows from 6.3 Proposition and (1.3.1) that 
$$
c(\psi_K) = d(K|T)+e(K|T)c(\psi_T) = p^r{-}2. 
$$ 
Let $\beta_K\in K^\times$ satisfy (1). The coset $\beta_KU^1_K$ is thereby uniquely determined, and $\ups_K(\beta_K) = -(1{+}p^r)$. 
\par 
Let $X_1(T)$ denote the group of tamely ramified characters of $T^\times$. Let $\eps\in X_1(T)$ and set $\eps_K = \eps\circ \N 
KT$. The map $\eps \mapsto \eps_K$ is an isomorphism $X_1(T)\to X_1(K)$, since $K/T$ is totally wildly ramified. The 
induction relation 
$$
\Ind_{K/T}\,\eps_K\otimes \chi = \eps\otimes \sigma_T, 
$$ 
implies the local constant relation 
$$ 
\frac{\ve(\eps\otimes \sigma_T,s,\psi_T)}{\ve(\sigma_T,s,\psi_T)} = \frac{\ve(\eps_K\otimes\chi,s,\psi_K)}{\ve(\chi,s,\psi_K)}\,. 
$$ 
By 5.5 Proposition and a calculation parallel to 2.2 Lemma, the element $\gamma_{\sigma_T}$ of 2.3 Lemma is $\det\alpha$ (modulo $U^1_T$). It follows that $\eps(\det\alpha) = \eps_K(\beta_K)$ for all $\eps\in X_1(T)$, and therefore $\det\alpha \equiv \N KT(\beta_K) \pmod{U^1_T}$. 
\par
Conversely, let $\beta_K'\in K^\times$ satisfy $\N KT\,\beta'_K \equiv \det\alpha \pmod{U^1_T}$. It follows that $\beta'_K \equiv \beta_K \pmod{U^1_K}$, so $\psi_K(\beta'_Ky) = \psi_K(\beta_Ky) = \chi(1{+}y)$, $y\in \frak p_K^2$. \qed 
\enddemo 
We re-formulate the proposition to avoid the choice of a Lagrangian. 
\proclaim{Corollary} 
There exists $\beta_E \in E^\times$ such that $\N ET\,\beta_E \in (\det\alpha)U^1_T$. For any such element $\beta_E$, we have 
$$ 
\xi_\sigma(1{+}z) = \psi_E(\beta_E^{p^r}z), \quad z\in \frak p_E^{1+p^r}, 
$$
where $\psi_E = \psi_F\circ \roman{Tr}_{E/F}$. 
\endproclaim 
\demo{Proof} 
The first assertion is immediate. If $\beta_K = \N EK\,\beta_E$, then $\det\alpha \equiv \N KT\,\beta_K \pmod{U^1_T}$ while $\beta_K \equiv \beta_E^{p^r} \pmod{U^1_E}$. By 6.4 Proposition, 
$$
\N EK(1{+}z) \equiv 1+\roman{Tr}_{E/K}(z) \pmod{\frak p_K^3}, \quad z\in \frak p_E^{1+p^r}, 
$$ 
so $\xi_\sigma(1{+}z) = \psi_K(\beta_K\roman{Tr}_{E/K}(z)) = \psi_E(\beta_K z) = \psi_E(\beta_E^{p^r}z)$, as required. \qed 
\enddemo  
\subhead 
6.6 
\endsubhead 
By definition, the character $\xi_\sigma$ is fixed under the action of $\vT = \Gal E{F'}$ and, by 6.3 Proposition, $\sw(\xi_\sigma) = 1{+}p^r$. We examine some implications of these properties. 
\par 
The soluble group $\vT$ has order $p^{2r}(1{+}p^r)$. The Hall Subgroup Theorem gives a subgroup $\vD$ of $\vT$ of order $1{+}p^r$, unique up to conjugation in $\vT$. If we set $H = E^\vD$, the extension $E/H$ is totally tamely ramified of degree $1{+}p^r$ and so $\vD$ is cyclic. 
\proclaim{Proposition} 
Let $\vD$ be a subgroup of $\vT$ of order $1{+}p^r$, and let $H = E^\vD$. 
\roster 
\item 
Let $\tau$ be a character of $E^\times$ fixed by $\vD$. There exists a character $\tau^H$ of $H^\times$ such that $\tau = \tau^H\circ \N EH$. Moreover, $\sw(\tau) = (1{+}p^r)\,\sw(\tau^H)$. 
\item 
The norm map $\N EH$ induces an isomorphism 
$$ 
\align 
U^{1+p^r}_E/U^{2+p^r}_E &\longrightarrow U^1_H/U^2_H, \\ 
1{+}x &\longmapsto 1+\roman{Tr}_{E/H}(x). 
\endalign  
$$
\item 
Let $\tau_1$, $\tau_2$ be $\vD$-fixed characters of $E^\times$ satisfying $\sw(\tau_1) = \sw(\tau_2) = 1{+}p^r$. If $\tau_1$, $\tau_2$ agree on $U^{1+p^r}_E$, then $\tau_2 = \phi\tau_1$, for a tamely ramified character $\phi$ of $E^\times$. 
\endroster 
\endproclaim 
\demo{Proof} 
Parts (1) and (2) are standard properties of cyclic, tamely ramified extensions. Part (3) follows from (1) and (2). \qed 
\enddemo  
\subhead 
6.7 
\endsubhead 
We interpolate an alternative description of the character $\xi_\sigma|_{U^1_E}$. Taking $\vD$ and $H = E^\vD$ as in 6.6, choose $\beta_H\in H$ so that 
$$
\N H{F'}\,\beta_H \equiv \det\alpha \pmod{U^1_{F'}}. 
\tag 6.7.1 
$$ 
The relation (6.7.1) determines $\beta_H$ modulo $U^1_H$, that is, modulo $U^{1+p^r}_E$. 
\par
Let $\psi_H = \psi_F\circ\roman{Tr}_{H/F}$, and define a character $\xi^H$ of $U^1_H/U^2_H$ by  
$$
\xi^H(1{+}x) = \psi_H(\beta_H^{p^r}x), \quad x\in \frak p_H. 
$$
Thus $\xi^H\circ \N EH$ is a character of $U^1_E$, fixed by $\vD$ and agreeing with $\xi_\sigma$ on $U^{1+p^r}_E$. That is, 
$$
\xi_\sigma(u) = \xi^H(\N EH(u)), \quad u\in U^1_E, 
\tag 6.7.2 
$$ 
by 6.6 Proposition and 6.5 Corollary. 
\subhead 
6.8 
\endsubhead 
We prove 6.1 Theorem. Property (1) is a direct consequence of the definition of $\xi_\sigma$, while (2) is given by 6.3 Proposition. We have just established (3) in 6.5 Corollary. The uniqueness assertion is given by 6.6 Proposition. \qed
\par 
We prove 6.1 Corollary. The equivalence of (1) and (2) is given by 6.5 Proposition, and that of (2) and (3) by 6.5 Corollary and 5.2 Corollary. \qed 
\remark{Remark} 
Theorem 6.1 gives $\xi_\sigma$ on $U^1_E$ while $\xi_\sigma^{p^r} = \det\sigma|_{\scr W_E} = \omega_\pi \circ \N EF$. Together, these relations determine $\xi_\sigma$ up to an unramified factor of order dividing $p^r$. One cannot isolate this factor via the standard technique of a local constant calculation. The character $\psi_E = \psi_F\circ \roman{Tr}_{E/F}$ has $c(\psi_E) = p^{2r}{-}2$ (by (1.3.2) and 6.3 Proposition). So, if $\chi$ is unramified and $\chi^{p^r} = 1$, then $\ve(\chi\xi_\sigma,s,\psi_E) = \ve(\xi_\sigma,s,\psi_E)$. 
\endremark 
\parskip=2pt plus.5pt minus.5pt
 
\enddocument 
\Refs 
\ref\no1
\by C.J. Bushnell 
\paper Hereditary orders, Gauss sums and supercuspidal representations of  $GL_N$ 
\jour J. reine angew. Math. \vol 375/376 \yr 1987 \pages184--210 
\endref 
\ref\no2
\by C.J. Bushnell and A. Fr\"ohlich 
\book Gauss sums and $p$-adic division algebras 
\bookinfo Lecture Notes in Math. \vol 987 \publ Springer \yr 1983 
\endref 
\ref\no3
\bysame 
\paper Non-abelian congruence Gauss sums and $p$-adic simple algebras 
\jour  Proc. London Math. Soc. (3) \vol 50  \yr 1985 \pages 207--264 
\endref 
\ref\no4 
\by C.J. Bushnell and G. Henniart 
\paper Local tame lifting for $GL(N)$ I: simple characters 
\jour Publ. Math. IHES \vol 83 \yr 1996 \pages 105--233 
\endref 
\ref\no5 
\bysame 
\paper Local tame lifting for $GL(n)$ II: wildly ramified supercuspidals 
\jour Ast\'erisque \vol 254 \yr 1999 \pages 1-105 
\endref 
\ref\no6
\bysame 
\paper Local tame lifting for $GL(n)$ IV: simple characters and base change 
\jour Proc. London Math. Soc. (3) \vol 87 \yr 2003 \pages 337--362 
\endref 
\ref\no7 
\bysame 
\paper The essentially tame local Langlands correspondence, I 
\jour J. Amer. Math. Soc. \vol 18 \yr 2005 \pages 685--710 
\endref 
\ref\no8 
\bysame 
\paper The essentially tame local Langlands correspondence, II: totally ramified representations 
\jour Compositio Mathematica \vol 141 \yr 2005\pages 979--1011 
\endref 
\ref\no9 
\bysame 
\book The local Langlands Conjecture for $GL(2)$ 
\bookinfo Grundlehren der mathematischen Wissenschaften {\bf 335} \yr 2006 \publ Springer 
\endref 
\ref\no10
\bysame 
\paper The essentially tame local Langlands correspondence, III: the general case 
\jour Proc. London Math. Soc. (3) \vol 101 \yr 2010 \pages 497--553 
\endref 
\ref\no11 
\bysame 
\paper Intertwining of simple characters in $\text{\rm GL}(n)$ 
\jour Internat. Math. Res. Notices; doi:10.1093/imrn/rns162 \yr 2012 
\endref 
\ref\no12 
\paper To an effective local Langlands Correspondence 
\jour Memoirs  Amer. Math. Soc., to appear. arXiv:1103.5316 \yr 2011 
\endref 
\ref\no13 
\by C.J. Bushnell, G.M. Henniart and P.C. Kutzko 
\paper Local Rankin-Selberg convolutions for $GL_n$: explicit conductor formula 
\jour J. Amer. Math. Soc. \vol 11 \yr 1998 \pages 703--730 
\endref 
\ref\no14 
\by C.J. Bushnell and P.C. Kutzko 
\book The admissible dual of $GL(N)$ via compact open subgroups 
\bookinfo Annals of Math. Studies \vol 129 \publ Princeton University Press \yr 1993 
\endref 
\ref\no15  
\bysame 
\paper The admissible dual of $SL(N)$ II 
\jour Proc. London Math. Soc. (3) \vol 68 \yr 1992 \pages 317--379 
\endref 
\ref\no16
\by P. Deligne and G. Henniart 
\paper Sur la variation, par torsion, des constantes locales d'equati\-ons fonctionnelles des fonctions $L$ 
\jour Invent. Math. \vol 64 \yr 1981 \pages 89--118 
\endref 
\ref\no17
\by R. Godement and H. Jacquet 
\book Zeta functions of simple algebras 
\bookinfo Lecture Notes in Math. \vol 260 \yr 1972 \publ Springer 
\endref 
\ref\no18 
\by B.H. Gross and M. Reeder 
\paper 
Arithmetic invariants of discrete Langlands parameters 
\jour Duke Math. J. \vol 154 \yr 2010 \pages 431--508 
\endref 
\ref\no19 
\by G. Henniart 
\paper R\'epresentations du groupe de Weil d'un corps local 
\jour L'Enseign. Math. \vol 26 \yr 1980 \pages 155--172 
\endref 
\ref\no20 
\bysame 
\book La conjecture de Langlands locale pour $GL(3)$ 
\bookinfo Memoire Soc. Math. France \vol 11/12 \yr 1984 \publ Gauthiers-Villars 
\endref 
\ref\no21  
\by H. Jacquet 
\paper Principal $L$-functions of the linear group  
\inbook Automorphic forms, representations and $L$-functions \eds A. Borel and 
W. Casselman \bookinfo Proc. Symposia Pure Math. \vol{XXXIII part 2} \publ Amer. Math. Soc. \yr 1979 \pages 63--86 
\endref 
\ref\no22 
\by H. Jacquet, I. I. Piatetskii-Shapiro and J. A. Shalika 
\paper Rankin-Selberg convolutions 
\jour Amer. J. Math. \vol 105 \yr 1983 \pages 367--483 
\endref 
\ref\no23 
\by H. Koch 
\paper Classification of the primitive representations of the Galois group of local fields 
\jour Invent. Math. \vol 40 \yr 1977 \pages 195--216 
\endref 
\ref\no24 
\by P.C. Kutzko 
\paper The Langlands conjecture for $GL_2$ of a local field 
\jour Ann. Math. \vol 112 \yr 1980 \pages 381--412 
\endref 
\ref\no25 
\by C. M\oe glin 
\paper Sur la correspondance de Langlands-Kazhdan 
\jour J. Math. Pures et Appl. (9) \vol 69 \yr 1990 \pages 175--226 
\endref 
\ref\no26 
\by 
\by M. Reeder and J.-K. Yu 
\paper Epipelagic representations and invariant theory 
\jour J. Amer. Math. Soc., to appear 
\endref 
\ref\no27 
\by J. F. Rigby 
\paper Primitive linear groups containing a normal nilpotent subgroup larger than the centre of the group 
\jour J. London Math. Soc. \vol 35 \yr 1960 \pages 389--400 
\endref 
\ref\no28 
\by J-P. Serre 
\book Corps locaux 
\publ Hermann \publaddr Paris \yr 1968 
\endref 
\ref\no29 
\by F. Shahidi 
\paper Fourier transforms of intertwining operators and Plancherel measures for $GL(n)$ \jour Amer. J. Math. \vol 106 \yr 1984 \pages 67--111
\endref 
\ref\no30 
\by J. Tate 
\paper Number theoretic background \inbook Automorphic forms, representations and $L$-func\-tions \eds A. Borel and W. Casselman \bookinfo Proc. Symposia 
Pure Math. \vol{XXXIII part 2}  \publ Amer. Math. Soc. \yr 1979 \pages 3--26 
\endref 
\endRefs
\bye